\documentclass{article}

\usepackage[utf8]{inputenc}

\usepackage{graphicx}
\usepackage{color}

\usepackage{pb-diagram}
\usepackage{hyperref}
\usepackage[leftcaption]{sidecap}
\sidecaptionvpos{figure}{t}
\usepackage{caption}

\usepackage{amsmath,amssymb,amsthm,amsfonts,amscd,euscript}

\title{A 15-vertex triangulation of the quaternionic projective plane}
\author{Denis Gorodkov\thanks{This work has been supported in part by the Moebius Contest Foundation for Young Scientists and by the Russian Science Foundation (project 14-50-00005).} \thanks{e-mail: \href{mailto:denis.gorod@mi.ras.ru}{denis.gorod at mi.ras.ru}}}
\date{}

 \newtheorem{proposition}{Proposition}
\newtheorem{lemma}{Lemma}
\newtheorem{theorem}{Theorem}
\newtheorem*{ntheorem}{Theorem}
\newtheorem{corollary}{Corollary}
\newtheorem{conjecture}{Conjecture}
\newtheoremstyle{neosn}{0.5\topsep}{0.5\topsep}{\rm}{}{\sc}{.}{ }{\thmname{#1}\thmnumber{ #2}\thmnote{ {\mdseries#3}}}
\theoremstyle{neosn}
\newtheorem{definition}{Definition}[section]
\newtheorem{remark}{Remark}[section]

\renewcommand\footnotemark{}

\newcommand{\link}{\mathrm{link}\,}
\newcommand{\Hom}{\mathrm{Hom}\,}

\begin{document}

\maketitle

\begin{abstract}
	In 1992, Brehm and K\"{u}hnel constructed a 8-dimensional simplicial complex $M^8_{15}$ with 15 vertices as a candidate to be a minimal triangulation of the quaternionic projective plane. They managed to prove that it is a manifold ``like a projective plane'' in the sense of Eells and Kuiper. However, it was not known until now if this complex is PL homeomorphic (or at least homeomorphic) to $\mathbb{H}P^2$. This problem was reduced to the computation of the first rational Pontryagin class of this combinatorial manifold. Realizing an algorithm due to Gaifullin, we compute the first Pontryagin class of $M^8_{15}$. As a result, we obtain that it is indeed a minimal triangulation of $\mathbb{H}P^2$.
\end{abstract}
	
\section{Introduction}
	
	A triangulation of a PL--manifold is a simplicial complex which is PL homeomorphic to this manifold. A~triangulation of a manifold is called \textit{(vertex-)minimal} if there are no triangulations of the same manifold with less vertices. The problem of finding minimal triangulations of a manifold is a classic problem in combinatorial topology; one can find a compilation of the most significant results on minimal triangulations in the survey by F. Lutz \cite{lutz}. The most interesting examples appear when the minimal triangulation has additional properties such as a non--trivial symmetry group. One of the well--known examples is the minimal triangulation of $\mathbb{R}P^2$ with $6$ vertices. It can be obtained by taking the quotient of the boundary of the icosahedron by the antipodal involution. In 1983 K\"uhnel and Banchoff constructed a simplicial complex named $\mathbb{C}P^2_9$ with $9$ vertices and proved that this complex is the minimal triangulation of the complex projective plane $\mathbb{C}P^2$. Besides, the symmetry group of this complex has order 54. Using similar ideas Brehm and K\"uhnel \cite{BK} constructed a 15-vertex simplicial complex $M^8_{15}$(as well as two other complexes $\widetilde{M}\vphantom{M}^8_{15}$ and $\widetilde{\widetilde{M}}\vphantom{M}^8_{15}$ that are PL homeomorphic to $M^8_{15}$) and conjectured that these complexes are minimal triangulations of the quaternionic projective plane $\mathbb{H}P^2$ where the PL structure on $\mathbb{H}P^2$ is induced by the canonical smooth structure.
	
	Brehm and K\"uhnel made an attempt to prove that the simplicial complex $M^8_{15}$ is PL homeomorphic to $\mathbb{H}P^2$, but they managed only to prove a weaker statement: $M^8_{15}$ is a manifold ``like a projective plane'', ie a manifold that admits a Morse function with exactly 3 critical points. Eells and Kuiper \cite{EK} examined this case in detail. In particular they showed that in dimension 8 such manifolds are distinguished by their Pontryagin numbers. Thus if we prove that Pontryagin numbers of the manifold $M^8_{15}$ coincide with Pontryagin numbers of $\mathbb{H}P^2$, this will imply that these manifolds are PL homeomorphic, i.e. $M^8_{15}$ is a triangulation of $\mathbb{H}P^2$. Moreover, Eells and Kuiper proved that for any 8-manifold ``like a projective plane'' its cohomology ring is isomorphic to the cohomology ring of $\mathbb{H}P^2$, i.e. $H^*(M^8_{15}, \mathbb{Z})=\mathbb{Z}[u]/(u^3)$, $\deg u = 4$. Hence, as an implication of Hirzebruch's formula for the signature of an 8--manifold, it is sufficient to compute the first rational Pontryagin class of $M^8_{15}$ to compute its Pontryagin numbers.
	
	As of the time Eells and Kuiper published their paper, there was no approach for computing the first Pontryagin class of a triangulated manifold that would be appropriate for explicit computations. Formulae that were known by that time (\cite{ggl,ggl2,ggl3,mcpherson,cheeger,gelmcph}) were not fully combinatorial, i.e. they did not give the posibility to compute the first Pontryagin class using only the combinatorial structure of the triangulation. Moreover, all these formulae require difficult and laborious computations. The only example of an explicit computation using one of these formulae -- the Gabrielov--Gelfand--Losik formula \cite{ggl, ggl3, ggl2} -- is the computation by Milin \cite{milin} of the first Pontryagin class of $\mathbb{C}P^2_9$.
	
	In 2004 Gaifullin \cite{G04} (cf. \cite{G08,bis}) constructed an explicit algorithm for computing the first rational Pontryagin class of a combinatorial manifold. \textit{A combinatorial manifold} of dimension $n$ is a simplicial complex $K$, such that the link of any vertex of $K$ is PL homeomotphic to the boundary of the $n$--dimensional simplex. Note that any PL triangulation of a PL manifold is a combinatorial manifold. This algorithm is fully combinatorial, i.e. the computation does not need any additional data except the combinatorial structure of the triangulation.

\section{Main results}
	
	\begin{theorem}
		\label{main}
		The first rational Pontryagin class $p_1(M^8_{15})$ is equal to $2u$ where $u$ is the image of one of the two generators of the group $H^4(M^8_{15},\mathbb{Z})\cong \mathbb{Z}$ under the natural embedding $H^4(M^8_{15},\mathbb{Z})\subset H^4(M^8_{15},\mathbb{Q})$.
	\end{theorem}
	
	\begin{remark}
		The results of Kervaire and Milnor \cite{kermil} imply that the groups of smooth structures on spheres modulo $h$-cobordism are trivial up to dimension $6$. An easy consequence from this fact is the following: unlike higher Pontryagin classes, the first integral Pontryagin class is a PL invariant and is well-defined for PL manifolds (cf. \cite{brumfiel}). Thus, our theorem can be reformulated in the following way: the first integral Pontryagin class $p_1(M^8_{15})$ is equal to $2u$ where $u$ is one of two generators of the group $H^4(M^8_{15},\mathbb{Z})$.
	\end{remark}
	
	\begin{corollary}\label{cor}
		Pontryagin numbers of the combinatorial manifold $M^8_{15}$ are equal to corresponding Pontryagin numbers of the quaternionic projective plane $\mathbb{H}P^2$.
	\end{corollary}
	
	\begin{proof}
		It follows from the classical Hirzebruch formula that the signature of a closed oriented manifold can be obtained as a linear combination of Pontryagin numbers of the manifold. For an 8-manifold it looks as follows:
		$$\sigma(X)=\dfrac{7p_2[X]-p_1^2[X]}{45}$$
		As the cohomology rings of $\mathbb{H}P^2$ and $M^8_{15}$ are isomorphic, we can choose the orientation on $M^8_{15}$ such that $\sigma(M^8_{15})=\sigma(\mathbb{H}P^2)=1$. It is well-known that $p_1^2[\mathbb{H}P^2] = 4$. Thus, $$p_1^2[M^8_{15}]=\langle(2u)\smile(2u), [M^8_{15}]\rangle=4=p_1^2[\mathbb{H}P^2]$$ Finally, $p_2[M^8_{15}]=p_2[\mathbb{H}P^2]=7$.
	\end{proof}
	
	It follows from Corollary \ref{cor} and from the results of \cite{BK, BK2, EK} that
	\begin{corollary}
		$M^8_{15}$, $\widetilde{M}\vphantom{M}^8_{15}$ and  $\widetilde{\widetilde{M}}\vphantom{M}^8_{15}$ are PL--homeomorphic to $\mathbb{H}P^2$ and are minimal triangulations of $\mathbb{H}P^2$.
	\end{corollary}
	
	This corollary will be accurately proved in the end of Section \ref{program}.
	
	The algorithm for computing the first Pontryagin class was implemented on a computer in the general case using the programming language \texttt{GAP}(\cite{GAP4}).
	
\section{Manifolds ``like a projective plane''}
\label{EK}
	
	The classical notion of a Morse function can be generalized for topological or combinatorial manifolds in the following way. (The author took this generalization from the article \cite{EK}. As of today, another non-equivalent definition of a Morse function on a combinatorial manifold is used. See \cite{forman} for the modern combinatorial Morse theory.)
	
	Consider one of the first assertions of classical Morse theory:
	\begin{proposition}
		Let $M^n$ be a smooth manifold. If $a \in M$ is a non-critical point of a Morse function $f \colon M \longrightarrow~\mathbb{R}$. Then there is a smooth $a$--centered coordinate system $\{ x^i \} $ such that $x^n = f(x)-f(a)$ in a neighbourhood of the point $a$. If $a$ is a critical point of $f$, then there is a smooth $a$--centered coordinate system $\{ x^i \} $ such that $$-\sum_{i=1}^k (x^i)^2 + \sum_{i=k+1}^n (x^i)^2 = f(x)-f(a)$$ in a neighbourhood of the point $a$.
	\end{proposition}
	
	This crucial statement can be taken as a definition of a Morse function in the smooth case. In the topological and combinatorial case one can use this approach, as we can give up the requirement of smoothness.
	
	In the case of a combinatorial manifold $K$ the function $f$ is meant as a function on the geometrical realization $|K|$ of the manifold.
	
	\begin{definition}
	\label{morse}
	A Morse function on a topological(respectively, combinatorial) manifold $X$ is a continuous (respectively, piecewise linear) function $f\colon X\longrightarrow \mathbb{R}$, such that in the neighbourhood of any point $a\in X$ there is a continuous $a$--centered coordinate system $\{ x^i \} $, such that one of the two conditions (1) and (2) is satisfied (respectively, (1) and (2')): 
		\begin{enumerate}
		\item $f(x)-f(a) = x^n$ ; such a point $a$ is called ordinary.
		\item $f(x)-f(a) = -\sum_{i=1}^k (x^i)^2 + \sum_{i=k+1}^n (x^i)^2$ ; such a point $a$ is called critical of index $k$ in the topologial case.
		\item[2'.] $f(x)-f(a)=-\max\{|x^1|,\ldots,|x^k|\}+\max\{|x^{k+1}|,\ldots,|x^n|\}$ ; such a point $a$ is called critical of index $k$ in the combinatorial case.
		\end{enumerate}
	\end{definition}
	
	\begin{remark}
	To compare, consider the modern definition of a combinatorial Morse function from \cite{forman}. Let $K$ be a simplicial complex, $S$ be the set of all simplexes of $K$ and $S_d$ be the set of simplexes of dimension $d$. A \textit{discrete} Morse function on $K$ is a function on $S$, such that for each $\sigma \in S_d$ $$\#\{\tau\in S_{d+1} | \tau\supsetneq\sigma\text{ and }f(\tau)\leqslant f(\sigma)\}\leqslant 1$$
	$$\#\{\nu\in S_{d-1}| \nu\supsetneq\sigma\text{ and }f(\nu)\geqslant f(\sigma)\}\leqslant 1$$
	Most of classical results for Morse theory stay true in the combinatorial case. This definition is more universal and more practical than the one we use. In the case of Eells--Kuiper's definition most of the results follow from constructing special deformations, and in the modern definition most results are first of all combinatorial.
	In our present work we will use only the definition from Eells and Kuiper's article \ref{morse}.
	\end{remark}
	
	Studying manifolds that allow Morse functions with few critical points is a natural problem. It is well-known that if there is a Morse function on a manifold with exactly two critical points, then the manifold is necessarily homeomorphic to a sphere. Eells and Kuiper showed that in the case of three critical points the results are quite more complicated.

	\begin{theorem}[Eells, Kuiper \cite{EK}] \label{EKth}
		Given a manifold $X$ with a Morse function $f \colon X \longrightarrow \mathbb{R}$ with precisely $3$ critical points.
		\begin{enumerate}
			\item
				Dimension and cohomology.
				\textup{The only possible dimensions of $X$ are $n=0,\,2,\,4,\,8,\,16$.}
				\textup{For $n=0$\, the space  $X$ consists of three points. For $n=2$\, the space is homeomorphic to the real projective plane  $\mathbb{R}P^2$\,. }
				
				\textup{The cohomology ring $H^*(X, \mathbb{Z})$ is isomorphic to the cohomology ring of the complex~($n=4$), quaternionic~($n=8$) or Cayley~($n\nolinebreak=\nolinebreak16$) projective plane, ie $H^*(X, \mathbb{Z})=\mathbb{Z}[u]/(u^3)$ where $\dim u = n/2$.}
				
			\item
	
				\textup{$X$ is a compactification of $\mathbb{R}^{n}$ by a sphere $S^{n/2}$\,.}
	
			\item
				Homotopy type.
				
				\textup{For $n=4$ only one homotopy type of $X~\mathbb{C}P^2$ is possible, for $n=8$ there are $6$ homotopy types, and for $n=16$ there are $60$ of them.}
				
			\item
				From the combinatorial point of view, \textup{there are infinitely many different possible manifolds ``like a projective plane'' in dimensions $n=8$ and $n=16$. They are classified by their Pontryagin numbers, ie if two such manifolds have equal Pontryagin numbers, then these manifolds are PL homeomorphic. Some of these manifolds do not admit a compatible smooth structure.}
			
			\item
				\textup{In the case of dimension $n=8$ let us present the results more precisely. The Pontryagin number $p_1^2$ of the manifold $X^8$ can take the following form $$p_1^2[X]=4 (2h-1)^2 ,$$ where $h$ is an integer parameter, that parameterizes all the $X^8$.  Combinatorial manifolds $X^8_h$ admit a compatible smooth structure if and only if $h\equiv 4j$ or $h\equiv 4j+1$ modulo $12$. Moreover, $X^8_{h_0}$ and $X^8_{h_1}$ belong to the same homotopy class if and only if either $h_0-\nolinebreak h_1\equiv\nolinebreak 0$, either $h_0+h_1\equiv\nolinebreak 1$ modulo $12$.}
					
		\end{enumerate}
		
	\end{theorem}
	
	This theorem makes the following definition natural.
	\begin{definition}
			A manifold \textit{``like a projective plane''} is a topological, smooth or combinatorial manifold, such that there exists, respectively, a continuous, smooth or piecewise linear Morse function with three critical points.
	\end{definition}

\section{Brehm--K\"{u}hnel complexes} \label{complex}

	K\"{u}hnel and Banchoff \cite{BK} constructed a special 9-vertex simplicial complex $\mathbb{C}P^2_9$. It has several remarkable properties:
	
	\begin{enumerate}
	\item
		Among all combinatorial 4-manifolds which are not homeomorphic to the sphere it has the least number of vertices.
	\item
		Any 3 vertices of $\mathbb{C}P^2_9$ span a simplex contained in the complex (this property is called 3-neighbourliness). Moreover, 5 vertices of $\mathbb{C}P^2_9$ span a simplex iff the remaining 4 vertices do not span a simplex.
	\item
		This complex is a vertex-minimal triangulation of $\mathbb{C}P^2$.
	\item
		$\mathbb{C}P^2_9$ has an automorphism group of order 54.
	\end{enumerate}

	In an attempt to find the minimal triangulation of the quaternionic projective plane Brehm and K\"{u}hnel \cite{BK} constructed three 15-vertex simplicial complexes $M^8_{15}$, $\widetilde{M}\vphantom{M}^8_{15}$ and $\widetilde{\widetilde{M}}\vphantom{M}^8_{15}$ with similar properties in the 8-dimensional case.
	
	\begin{enumerate}
	\item
		Among all combinatorial 8-manifolds which are not homeomorphic to the sphere they have the least number of vertices.
	\item
		Any 5 vertices of any of the complexes $M^8_{15}$, $\widetilde{M}\vphantom{M}^8_{15}$, $\widetilde{\widetilde{M}}\vphantom{M}^8_{15}$ span a simplex(these complexes are 5-neighbourly). Moreover, 9 vertices span a simplex iff the remaining 6 vertices do not span a 5-simplex.
	\item
		The automorphism groups of $M^8_{15}$, $\widetilde{M}\vphantom{M}^8_{15}$ and $\widetilde{\widetilde{M}}\vphantom{M}^8_{15}$ are isomorphic to $A_5$, $A_4$ and $S_3$ respectively.
	
	\end{enumerate}
	
	The construction of the complexes is based on explicit descriptions of some group actions on the set of vertices. The actions will be given as subgroups of the permutation group on 15 elements $S_{15}$.
	
	Consider the following permutations:\\ 
	$$P = (1~2~3~4~5)(6~7~8~9~10)(11~12~13~14~15)$$
	$$T = (3~10)(4~14)(5~8)(6~11)(7~12)(13~15)$$
	$$U = (1~6~11)(2~7~12)(3~8~13)(4~9~14)(5~10~15)$$
	
	We will also need
	$$S = (1~6~11)(2~15~14)(3~13~8)(4~7~5)(9~12~10) = P^{-1}TP^{-2}TP^{-2}$$
	$$R = (2~5)(3~4)(7~10)(8~9)(12~15)(13~14) = S^{-1}P^2SP^{-1}S$$
	
	Then define $\mathcal{G}_2 = \langle P,T\rangle$, $\mathcal{G}_3 = \langle P,T,U\rangle$, $\mathcal{G}_1 = \langle P,S\rangle$, $\mathcal{G}_0 = \langle R,S\rangle$, $\widetilde{\mathcal{G}}_0 = \langle PRP^{-1},S\rangle$.
	
	We have the following natural injective homomorphisms:
	\[	\begin{diagram}
			\node{\mathcal{G}_0} \arrow{ese}\\  \node[3]{\mathcal{G}_1}\arrow{e}
			\node{\mathcal{G}_2}\arrow{e}
			\node{\mathcal{G}_3}
			\\
			\node{\widetilde{\mathcal{G}}_0} \arrow{ene,T}
			\end{diagram}
	\]
	
	The group $\mathcal{G}_1\cong A_5$ will be the automorphism group of $M^8_{15}$, and groups $\mathcal{G}_0$ and $\widetilde{\mathcal{G}}_0$ will be automorphism groups of $\widetilde{M}\vphantom{M}^8_{15}$ and $\widetilde{\widetilde{M}}\vphantom{M}^8_{15}$, respectively.
	
	The complexes $M^8_{15}$, $\widetilde{M}\vphantom{M}^8_{15}$ and $\widetilde{\widetilde{M}}\vphantom{M}^8_{15}$ consist of two parts: they have a common part $\mathcal{K}_0$ consisting of 415 simplices of maximal dimension, that is described as the union of orbits of 12 explicitely given simplexes under the action of $\mathcal{G}_1$: 	

	\begin{align*}
	A &= \{1, 2, 3, 6, 8, 11, 13, 14, 15\}\\
	B &= \{1, 3, 6, 8, 9, 10, 11, 12, 13\}\\
	C &= \{1, 2, 6, 9, 10, 11, 12, 14, 15\}\\
	D &= \{1, 2, 3, 4, 7, 9, 12, 14, 15\}\\
	E &= \{1, 2, 4, 7, 9, 10, 12, 13, 14\}\\
	F &= \{1, 2, 6, 8, 9, 10, 11, 14, 15\}\\
	G &= \{1, 2, 3, 4, 5, 6, 9, 11, 13\}\\
	H &= \{1, 3, 5, 6, 8, 9, 10, 11, 12\}\\
	I &= \{1, 3, 5, 6, 7, 8, 9, 10, 11\}\\
	J &= \{1, 2, 3, 4, 5, 7, 10, 12, 15\}\\
	K &= \{1, 2, 3, 7, 8, 10, 12, 13, 14\}\\
	M &= \{2, 5, 6, 7, 8, 9, 10, 13, 14\}\\
	\end{align*}
	
	To define the remaining 75\, 8-simplices for each of the complexes, consider the simplexes
	\begin{align*}
	L_{(1)} &= \{3, 4, 6, 7, 11, 12, 13, 14, 15\}\\
	N_{(1)} &= \{3, 4, 6, 7, 10, 12, 13, 14, 15\}\\
	\end{align*}
	and take their images under powers of the permutation $P$.
	
	$$L_{(n)}=P^{n-1}L_{(1)},~~~ N_{(n)}:=P^{n-1}N_{(1)}$$
	$$\widetilde{L}_{(n)}=P^{n-1}TL_{(1)},~~~ \widetilde{N}_{(n)}:=P^{n-1}TN_{(1)}$$
	
	Finally, denote $\mathcal{L}_n=L_{(n)}\cup N_{(n)}$ and $\mathcal{\widetilde{L}}_n=\widetilde{L}_{(n)}\cup \widetilde{N}_{(n)}$.
	
	$$\mathcal{K}_1=\mathcal{L}_1\cup \mathcal{L}_2\cup\mathcal{L}_3\cup\mathcal{L}_4\cup\mathcal{L}_5;$$
	$$\mathcal{\widetilde{K}}_1=\mathcal{\widetilde{L}}_1\cup \mathcal{L}_2\cup\mathcal{L}_3\cup\mathcal{L}_4\cup\mathcal{L}_5;$$
	$$\mathcal{\widetilde{\widetilde{K}}}_1=\mathcal{\widetilde{L}}_1\cup \mathcal{L}_2\cup\mathcal{\widetilde{L}}_3\cup\mathcal{L}_4\cup\mathcal{L}_5.$$

	Then the required complexes can be written as
	$$M^8_{15}=\mathcal{K}_0\cup\mathcal{K}_1;~~~\widetilde{M}\vphantom{M}^8_{15}=\mathcal{K}_0\cup\mathcal{\widetilde{K}}_1;~~~\widetilde{\widetilde{M}}\vphantom{M}^8_{15}=\mathcal{K}_0\cup\mathcal{\widetilde{\widetilde{K}}}_1.$$
	
	These three complexes are combinatorial manifolds and are PL homeomorphic to each other.
	
	The dimensions $n=2,4,8,16$ also appear precisely in the work \cite{BK2}. The authors consider all combinatorial manifolds and study the constraints that the number of vertices of a manifold apply on its dimension.
		
	\begin{theorem}[Brehm, K\"{u}hnel \cite{BK2}]
	\label{few_vert}
		Let $M^d$ be a combinatorial manifold with $n$ vertices. Then if $n<~\left\lceil\dfrac{3d}{2}\right\rceil + 3$, then $M$ is PL homeomorphic to a sphere, and if $n=\dfrac{3d}{2} + 3$, then either $M$ is PL homeomorphic to a sphere, either $M$ is a manifold ``like a projective plane''.
	\end{theorem}
		
	\begin{corollary}[\cite{BK}]
	$M^8_{15}$ is a manifold ``like a projective plane''.
	\end{corollary}
	\begin{proof}
	The complex $M^8_{15}$ is not homeomorphic to the sphere as its cohomology ring is not isomorphic to the cohomology ring of the sphere. So we can use Theorem \ref{few_vert}~, thus $M^8_{15}$ is a manifold ``like a projective plane''.
	\end{proof}
	
	In their article Brehm and K\"{u}hnel conjecture the following:
	\begin{conjecture}[\cite{BK}]
	$M_{15}^8$ is PL homeomorphic to $\mathbb{H}P^2$.
	\end{conjecture}
	
	Our main goal is to prove this conjecture.

	It follows from Theorems \ref{few_vert} and \ref{EKth}, that
	\begin{proposition}
	\label{prop}
		If $p_1^2\left[M^8_{15}\right]=p_1^2\left[\mathbb{H}P^2\right]$, than $M^8_{15}$ is PL homeomorphic to $\mathbb{H}P^2$ and is a minimal triangulation of it.
	\end{proposition}
	
	Let us now describe the algorithm of computing the first Pontryagin class $p_1(M^8_{15})$.

\section{Gaifullin's algorithm of computing the first Pontryagin class}
	
	The results of this section come from Gaifullin's article \cite{G04}.
	
	Denote by $\mathcal{T}_n$ the abelian group, generated by all isomorphism classes of oriented combinatorial $(n-1)$-spheres (ie $\left<L_1\right>=\left<L_2\right>$ if $L_1 \cong L_2$,  $\left<L\right>$ is the notation for the equivalence class of the sphere $L$) with relations $\left<-L\right>=-\left<L\right>$ where $-L$ is the notation for the sphere $L$ with reversed orientation.
	
	Let $f\in \Hom(\mathcal{T}_n,\mathbb{Q})$, and let $K^m$ be an oriented combinatorial manifold. Then denote $$f_\sharp(K)=\sum_{\sigma \in K\,, \dim \sigma = m-n} f(\left<\link\sigma\right>)\sigma.$$
		
	\begin{definition}
		A function $f\colon \mathcal{T}_n\longrightarrow\mathbb{Q}$ is a \textit{local formula} for a homogeneous polynomial $F \in \mathbb{Q}[p_1,p_2,\ldots]$ if for every $K$ the chain $f_\sharp(K)$ is a cycle, such that its homology class is dual to the class $F(p_1(K),p_2(K),\ldots)$.
	\end{definition}
	
	That is, the coefficient of a simplex $\sigma$ depends only on the combinatorial type of its ``neighbourhood'' -- the link.
	
	Our aim is the formula for the first Pontryagin class $f\colon \mathcal{T}_4\longrightarrow\mathbb{Q}$.

	\subsection{Bistellar moves}

		Let $K$ be a combinatorial manifold.
		
		\begin{definition}
			
			Let $\tau$ be a simplex, such that $\tau \notin K$, but all its faces lie in $K$ (we will call such a simplex \textit{empty}). Let also $sigma \in K$ be a simplex, such that $\sigma * \partial \tau$ is a complex of full dimension in $K$. Then \textit{a flip(also called a bistellar move or a Pachner move)} is a transformation of $K$ that replaces the subcomplex $\sigma * \partial \tau$ by $\tau * \partial \sigma$. We will also denote $\beta = \beta_{K, \sigma}$ and call $\beta$ the bistellar move, \textit{associated} with $\sigma$.
		
		\end{definition}
		
\begin{figure}[h]
			\centering
				\includegraphics[scale=0.6]{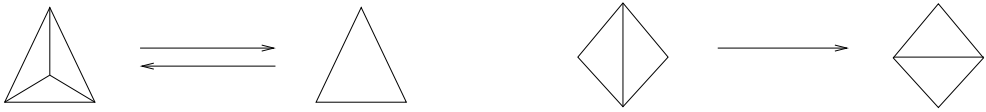}
				\caption{Moves in dimension 2}
				\label{flip}
		\end{figure}

		The bistellar moves in dimension 2 are shown on Fig.\ref{flip}.
			
		\begin{theorem}[Pachner, \cite{flip}]
			Let $K_1$ and $K_2$ be PL homeomorphic combinatorial manifolds. Then the manifold $K_1$ can be transformed in the manifold $K_2$ by a finite composition of bistellar moves and isomorphisms.
		\end{theorem}
		
		Thus, any two combinatorial spheres of the same dimension are connected by a sequence of bistellar moves. Then it is sufficient to show, how does the value $f(\left<L\right>)$ change(where $L$ is a combinatorial 3-sphere) under bistellar moves.
		
		Let $\beta_{K_1,\sigma} \colon K_1 \longrightarrow K_2$ and $\beta_{K'_1,\tau} \colon K'_1 \longrightarrow K'_2$ be bistellar moves. They are named \textit{equivalent} if there are isomorphisms $f\colon K_1 \longrightarrow K'_1$ and $f'\colon K_2 \longrightarrow K'_2$, such that $f(\sigma)=\tau$. If a bistellar move is equivalent to its inverse, we will call it \textit{inessential}, otherwise we will call the move essential.
						
		\subsection{The graph $\Gamma_2$}
		\label{Gamma2}
		
		Let us define a new construction -- the graph $\Gamma_n$. The vertices of this graph are oriented combinatorial spheres of dimension $n$. Two vertices $L_1$ and $L_2$ are connected with an edge if there is an essential bistellar move $\beta\colon L_1 \longrightarrow L_2$. If there are several non-equivalent bistellar moves between two vertices $L_1$ and $L_2$, then there are as many edges connecting $L_1$ and $L_2$ as equivalence classes of bistellar moves $\beta \colon L_1 \longrightarrow L_2$.
		
		Now consider, how should the value of the first Pontryagin class formula change is we transform a combinatorial sphere using a bistellar move. Let $\beta \colon L_1\longrightarrow L_2$ be a bistellar move, where $L_1$ and $L_2$ are combinatorial 3-spheres.
		
		Let $v$ be a vertex of $L_1$. Then we can consider the transformation induced by $\beta$ on the 2-sphere $\link_{L_1}v$. It is easy to show that this transformation is a bistellar move between 2-spheres $\link_{L_1}v\longrightarrow\link_{L_2}v$. Denote such a move by $\beta_v$.
		
		Gaifullin \cite{G04} constructed a special cohomology class $c\in H^1(\Gamma_2, \mathbb{Q})$ and proved the following theorem(for the explicit construction of the class $c$ see the table lower):
		
		\begin{theorem}[Gaifullin, \cite{G04}]
		\label{gaif}
			If $f\colon\mathcal{T}_4\longrightarrow\mathbb{Q}$ is a local formula for the first Pontryagin class, then for each bistellar move $\beta\colon L_1\longrightarrow L_2$ the following relation holds true:
			$$f(L_2)-f(L_1) = \sum_{v} h(\beta_v)\,,$$
			where $h \in \mathrm{C}^1(\Gamma_2,\mathbb{Q})$ is a cocycle of the graph $\Gamma_2$, representing the cohomology class $c$. For each cocycle $h \in C^1(\Gamma_2;\mathbb{Q})$, representing the cohomology class $c$, it is possible to explicitly indicate the function $f \in \Hom(\mathcal{T}_4,\mathbb{Q})$, which is a local formula for the class $p_1$.
		\end{theorem}
		
		To describe explicitly any local formula one should choose a precise representative $h$ of the cohomology class $c$. In \cite{G04} it is done in the following way. Let us choose for each vertex $L$ of the graph $\Gamma_2$ a chain $\xi$, such that $\partial\xi_{\{L\}} = \{L\}-\{\partial\Delta^3\}$. Consider all bistellar moves $\beta_1, \beta_2, \ldots \beta_r$ that lower the complexity (definition in the beginning of Section \ref{alg}) of the combinatorial sphere $L$, where for each $i$  $\beta_i \colon L \longrightarrow L_i$. Then assume $$\xi_{\{L\}} = \sum_{i=1}^r (\xi_{\{L_j\}} - \{\beta_j\}).$$ The desired cocycle is written by the formula $$h(\{\beta\})=\left<c, \{\beta\} + \{\xi_{L_1}\} - \{\xi_{L_1}\}\right>.$$
		Remark that this choice of the cocycle keeps the formula local, as it depends only on the combiantorial type of $L$. To describe explicitly the value of the cohomology class $c$ on a cycle in the graph $\Gamma_2$ we have to choose a set of linear generators among all cycles in $\Gamma_2$.
		
		\newpage
		
\begin{figure}[!hb]
			\centering
			\includegraphics[scale=0.55]{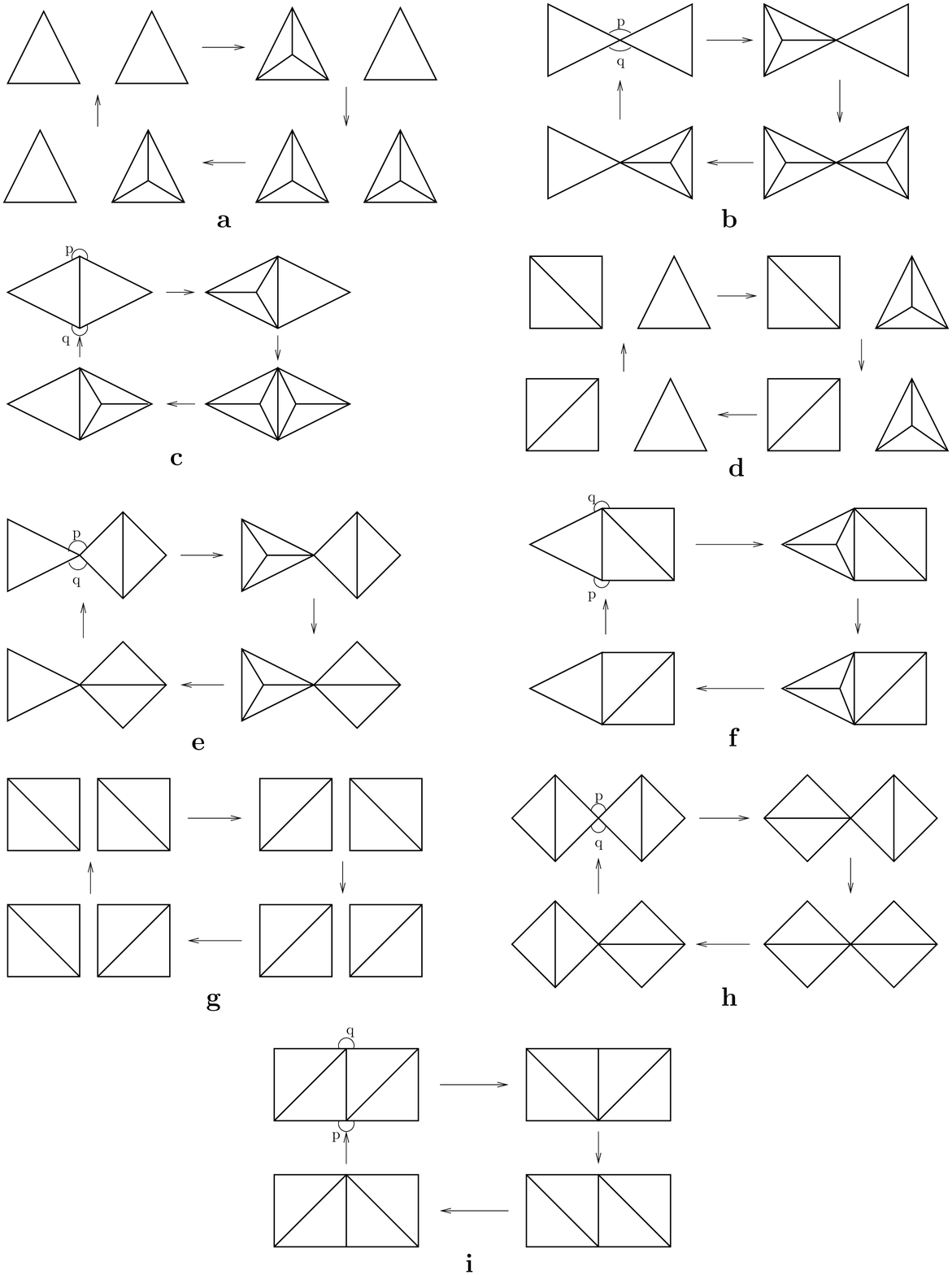}
			\caption{Elementary cycles of the first type}
			\label{pic-elem}
		\end{figure}

\begin{figure}[!hb]
			\centering
			\includegraphics[scale=0.4]{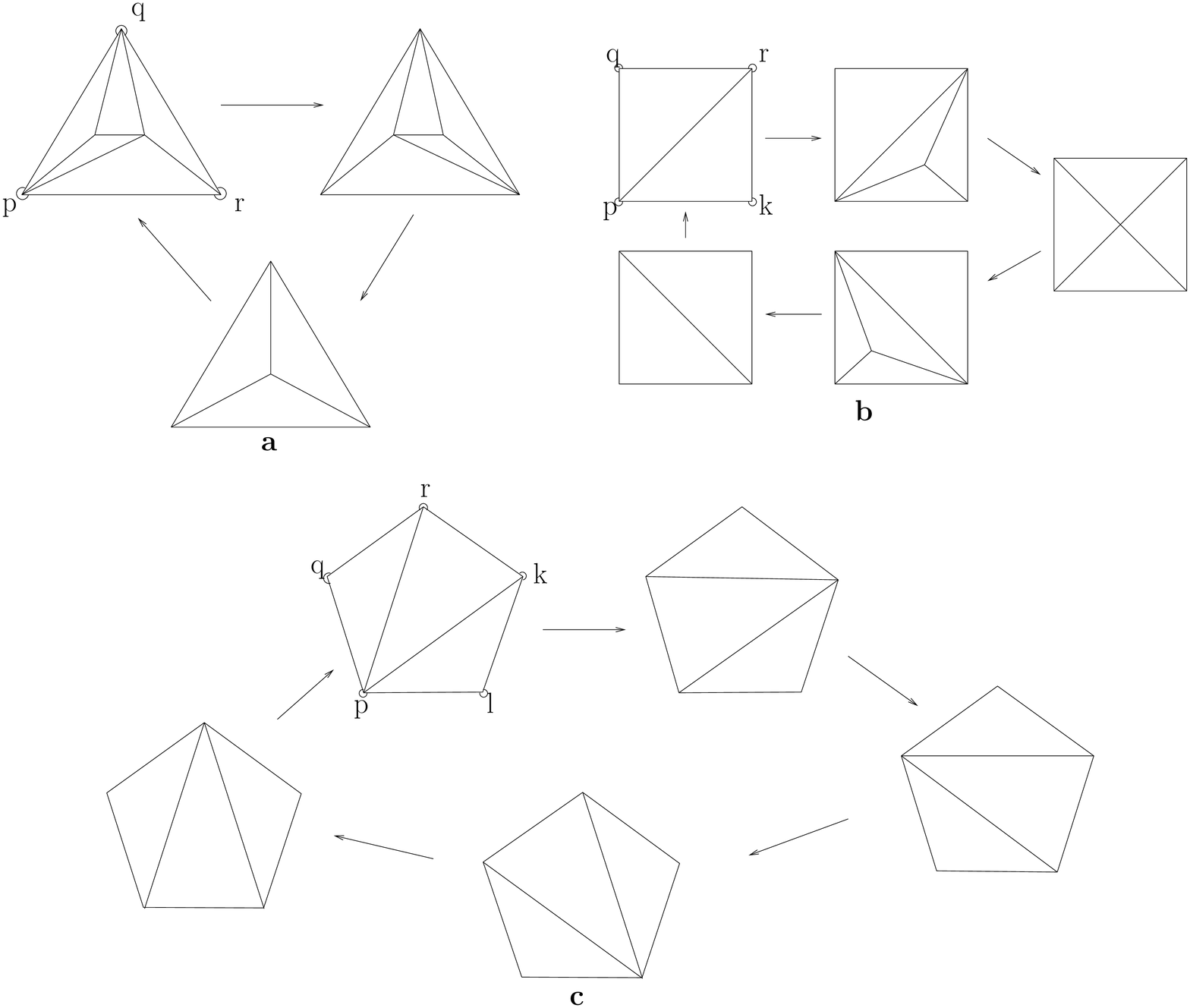}
			\caption{Elementary cycles of the second type}
			\label{pic-elem-2}
		\end{figure}

		So, a cycle in $\Gamma_2$ is given, ie a closed sequence of bistellar moves. We will call \textit{elementary cycles} of the first and second type some special cycles in the graph $\Gamma_2$. Cycles of the first and second type are shown in Fig.\ref{pic-elem} and Fig.\ref{pic-elem-2} respectively.

		The values of the cohomology class $c$ on elementary cycles were constructed in \cite{G04}, they depend on the number of triangles neighbouring the vertices whose link change after the bistellar moves. Numbers $p$, $q$, $r$, $k$, $l$ on Figures \ref{pic-elem} and \ref{pic-elem-2} denote the number of triangles inside the selected angles that contain the corresponding vertex.
		
		Consider two functions:
		$$\rho(p,q) = \dfrac{q-p}{(p+q+2)(p+q+3)(p+q+4)}$$
		$$\omega(p) = \dfrac{1}{(p+2)(p+3)}$$
		
		Then the value of the cohomology class $c$ on elementary cycles is given by the following table: \\		
		\begin{center}
		\begin{tabular}{| l | c |}
			\hline
			Type 1: a, d, g & 0 \\ \hline
			Type 1: b, e, h & $\rho(p,q)$ \\ \hline
			Type 1: c, i & $\rho(0,q) - \rho(0,p)$ \\ \hline
			Type 1: f & $\rho(0,q) + \rho(0,p)$ \\ \hline
			Type 2: a & $\omega(p) - \omega(q) + \omega(r) - \frac{1}{12}$ \\ \hline
			Type 2: b & $\omega(p) - \omega(q) - \omega(r) + \omega(k)$ \\ \hline
			Type 2: c & $\omega(p) + \omega(q) + \omega(r) + \omega(k) + \omega(l) - \frac{1}{12}$ \\ \hline
		\end{tabular}
		\end{center}
			
	\newpage
	
	Theorem \ref{gaif} provides the following algorithm for computig the first Pontryagin class.
			
			Consider a simplicial complex $K$. Choose its orientation.
			
			The algorithm of computing the first rational Pontryagin class consists of the following steps:
			\begin{enumerate}
				\item For each $(n-4)$-simplex of the complex $K$, find a sequence of bistellar moves that transform the link of this $(n-4)$-simplex into the boundary of a simplex.
						
				\item For each vertex $v$ of the link of each $(n-4)$-simplex consider $\link_{\link\sigma}(v)$. Then all obtained complexes are combinatorial 2-spheres. Induce the sequences of flips on these complexes as on subcomplexes of $\link\sigma$.
							
				\item For each obtained chain of flips, that reduce the combinatorial 2-sphere to the boundary of the 3-simplex, close the chain into a cycle into the complex $\Gamma_2$ in any way that depends only on the combinatorial type of the initial sphere.
				
				\item The resulting cycles are cycles in the graph $\Gamma_2$. Decompose this cycles in linear combinations of elementary cycles.
							
				\item Count the investment of each elementary cycle, recieve for each $\sigma$ the number $f\left( \left\langle link\,\sigma\right\rangle \right) $ and construct the cycle $$f_\sharp(K) = \sum_{\sigma \in K, \dim \sigma = n-4} f\left( \left\langle link\,\sigma\right\rangle \right)\, \sigma,$$ representing the homology element that is dual to the first Pontryagin class.
				
			\end{enumerate}
			
			The only remaining unexplained step is the decomposition of cycles in the graph $\Gamma_2$ into linear combinations of elementary cycles.			

	\subsection{Decomposition of cycles in the graph $\Gamma_2$ into linear combinations of elementary cycles}
	\label{alg}
		This algorithm was found by Gaifullin \cite{bis}, but some subcases were missed out. In the present article we eliminate the gap, thereby the algorithm of decomposition is now complete.
		
		We will often use the following notation. Suppose that $\sigma_1$ and $\sigma_2$ are simplexes in $L$, such that flips, associated with them are defined, and there is no simplex in $L$, containing both $\sigma_1$ and $\sigma_2$. Then denote by $\gamma(L,\sigma_1,\sigma_2)$ the following cycle:
				
				$$
				\begin{CD}
					L @>\beta_{\sigma_1}>> L_1\\
				  	@AA\beta^{-1}_{\sigma_2}A	@VV\beta_{\sigma_2}V\\
					L_3 @<\beta^{-1}_{\sigma_1}<< L_2
				\end{CD}
				$$
		
		We will say that the simplex \textit{participates} in the bistellar move, if the link of this simplex changes under the induced transformation.
		
		\begin{definition}
		The \textit{degree} of a vertex $v$ of a simplicial complex $K$ is the number of edges adjacent to this vertex.
		\end{definition}
		
		Let us introduce the notion of \textit{complexity} of a vertex of the graph $\Gamma_2$ as a combinatorial 2-sphere $L$ with $k$ vertices. 
		
		\[
			a(L)=
			\begin{cases}
				k, & \text{if $L$ contains at least one vertex of degree 3;} \\
				k+\frac{1}{3}, & \text{if $L$ contains a vertex of degree 4, but does not contain vertices of degree 3;} \\
				k+\frac{2}{3}, & \text{if $L$ does not contain vertices of degree 3 and 4.}
		    \end{cases}
		\]
		
		Now define the \textit{complexity} for edges of the graph $\Gamma_2$ (ie bistellar moves) $\beta\colon L_1 \rightarrow L_2$.
		
		\[
			a(\beta)=
			\begin{cases}
				\max(a(L_1),a(L_2)), & \text{if $a(L_1)\neq a(L_2)$\,;} \\
				a(L_1) + \frac{1}{6}, & \text{if $a(L_1)=a(L_2)$\,.}
			\end{cases}
		\]
		
		Then the complexity of any combinatorial sphere $a(L) \in \frac{1}{3}\mathbb{Z}_{\geqslant 0}$, and the complexity of any bistellar move $a(\beta) \in \frac{1}{6}\mathbb{Z}_{\geqslant 0}$.
		
		Denote the subgraph of the graph $\Gamma_2$, consisting of all vertices and edges with complexity not exceding $a$, by $\Gamma_2^a$ . Then if all cycles lying in the graph $\Gamma_2^a$ will be represented as a sum of a cycle from $\Gamma_2^{a-\frac{1}{6}}$ and elementary cycles, then by induction we will be able to represent all the cycle as a linear combination of elementary cycles. The base of the induction is the empty cycle for complexity of the bistellar moves that is equal to $4\frac{1}{6}$.
		
		On each step we will consider the least possible $a$ for a cycle.
		
		Let $a=k+\frac{b}{6}$. Then it is sufficient to prove the induction step for each of $b = 0, 1, 2, 3, 4, 5$. Consider separately the cases of even and odd $b$. Each of the cases is illustrated on the right by a corresponding image. Bistellar moves drawn on the picture are from the initial cycle and from the cycle with a smaller complexity as well as some auxiliary moves, in some cases the elementary cycles used are also denoted and shown on the picture.
		
		\paragraph{The case of odd $b$\,.}
			
			If $b$ is odd the transformations with the biggest complexity in the cycle $$\beta \colon L_1 \longrightarrow L_2,$$ join two combinatorial spheres with the same complexity $k + \frac{b-1}{6}$. If we are able to decompose each of these transformations into a linear combination of less complex bistellar moves and elementary cycles, then all the cycle can be represented as a sum of a cycle from $\Gamma_2^{a-\frac{1}{6}}$ and elementary cycles. We will process this step for each of the moves with the biggest complexity separately. Remark that the number of vertices in $L_1$ and $L_2$ coincide, hence $\beta$ is associated with an edge(not with a vertex). We denote this edge by $\sigma$.
		
			\subparagraph{$b=1$.}
				If $b=1$ both combinatorial spheres $L_1$ and $L_2$ contain vertices of degree 3.
				
				\begin{SCfigure}[50][!hb]
					\captionsetup{labelformat=empty}
					\includegraphics[scale=0.5]{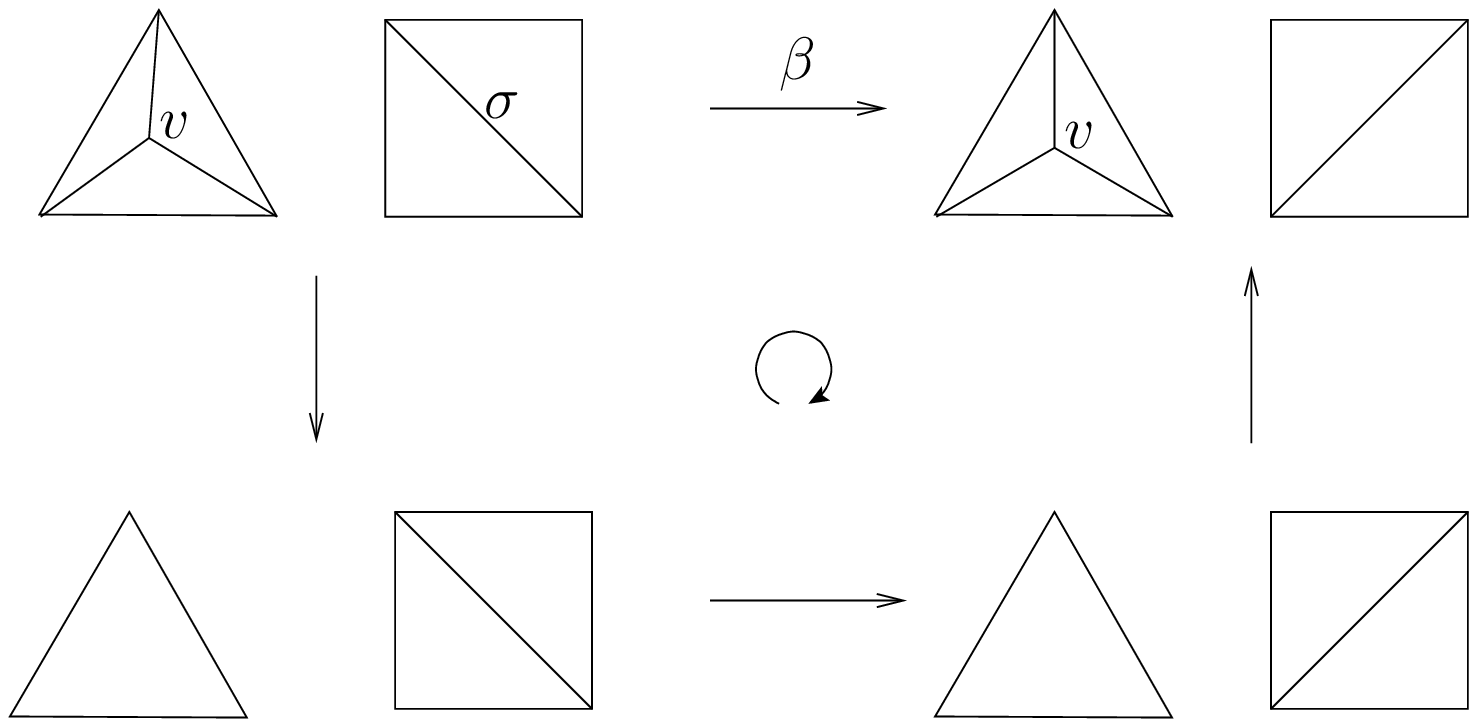}
					\caption{Suppose that there is a vertex $v$, such that its degree is equal to 3 in both spheres $L_1$ and $L_2$. In this case the cycle $\gamma(L_1,\sigma,v)$ is well-defined. The support of the chain $\beta-\gamma(L_1,\sigma,v)$ lies in the graph~$\Gamma_2^k$. }
				\end{SCfigure}
				
				\begin{SCfigure}[50][!hb]
					\captionsetup{labelformat=empty}
					\includegraphics[scale=0.45]{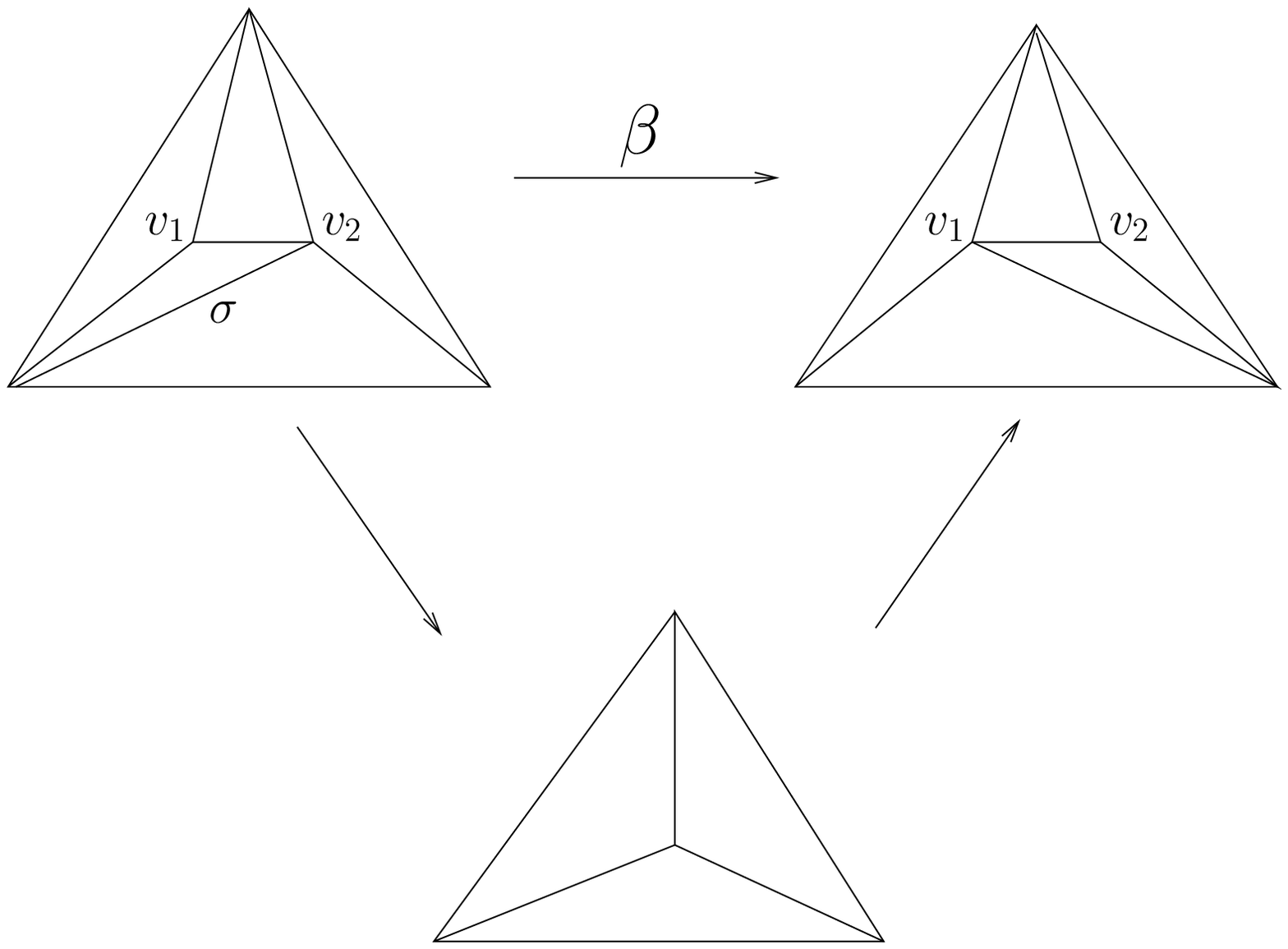}
					\caption{If there is no such vertex(with degree 3 in both spheres), then there are such vertices $v_1$ and $v_2$, that the degree of $v_1$ in $L_1$ and $L_2$ is equal to 3 and 4 respectively, and the degree of $v_2$ is equal to 4 and 3 respectively. Moreover, these vertices are joined by an edge, because the degree of both of them changed (one of them increasing, the other one decreasing) after one bistellar move $\beta$ associated with the edge $\sigma$. Then the cycle $\gamma(L_1,\sigma,v)$ is not defined. In this case a cycle $\delta$ of the second type is defined, and the support of the chain $\beta-\delta$ lies in the graph ~$\Gamma_2^k$.}
				\end{SCfigure}

			\subparagraph{$b=3$.}
							
				If $b=3$ both combinatorial spheres $L_1$ and $L_2$ do not contain vertices of degree 3, but contain a vertex of degree 4. Let us split in two cases.
				\newpage
				\begin{SCfigure}[50][!hb]
					\captionsetup{labelformat=empty}
					\includegraphics[scale=0.5]{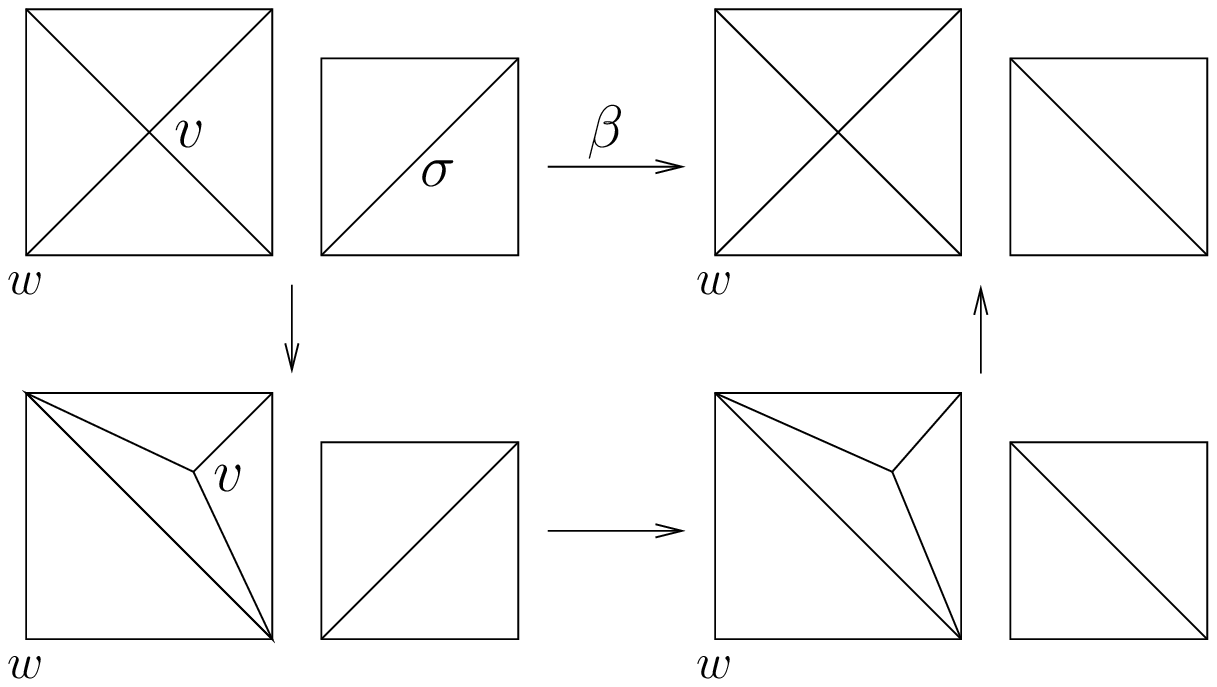}
					\caption{\textbf{1. Both combinatorial spheres $L_1$ and $L_2$ contain a common vertex of degree 4.} We will denote it by $v$. Then the vertex $v$ does not participate in the move $\beta$. Consider the tetragon $\link v$. $L_1$ and $L_2$ can not contain both diagonals of this tetragon. The move $\beta$ can be associated with the diagonal, but it can not replace one diagonal of the tetragon $\link v$ with the other (if it is the case, then $L_1$ has 5 vertices and the move $\beta$ is inessential). Thus, there is a diagonal of $\link v$ that is not contained in both $L_1$ and $L_2$.}
				\end{SCfigure}
				
				Denote one of the vertices of $\link v$, not belonging to this diagonal, by $w$. In this case the cycle $\gamma(L_1,\sigma,vw)$ is defined, and the support of the chain $\beta-\gamma(L_1,\sigma,vw)$ lies in $\Gamma_2^{k+\frac{1}{3}}$.

				\begin{SCfigure}[50][!hb]
					\captionsetup{labelformat=empty}
					\includegraphics[scale=0.6]{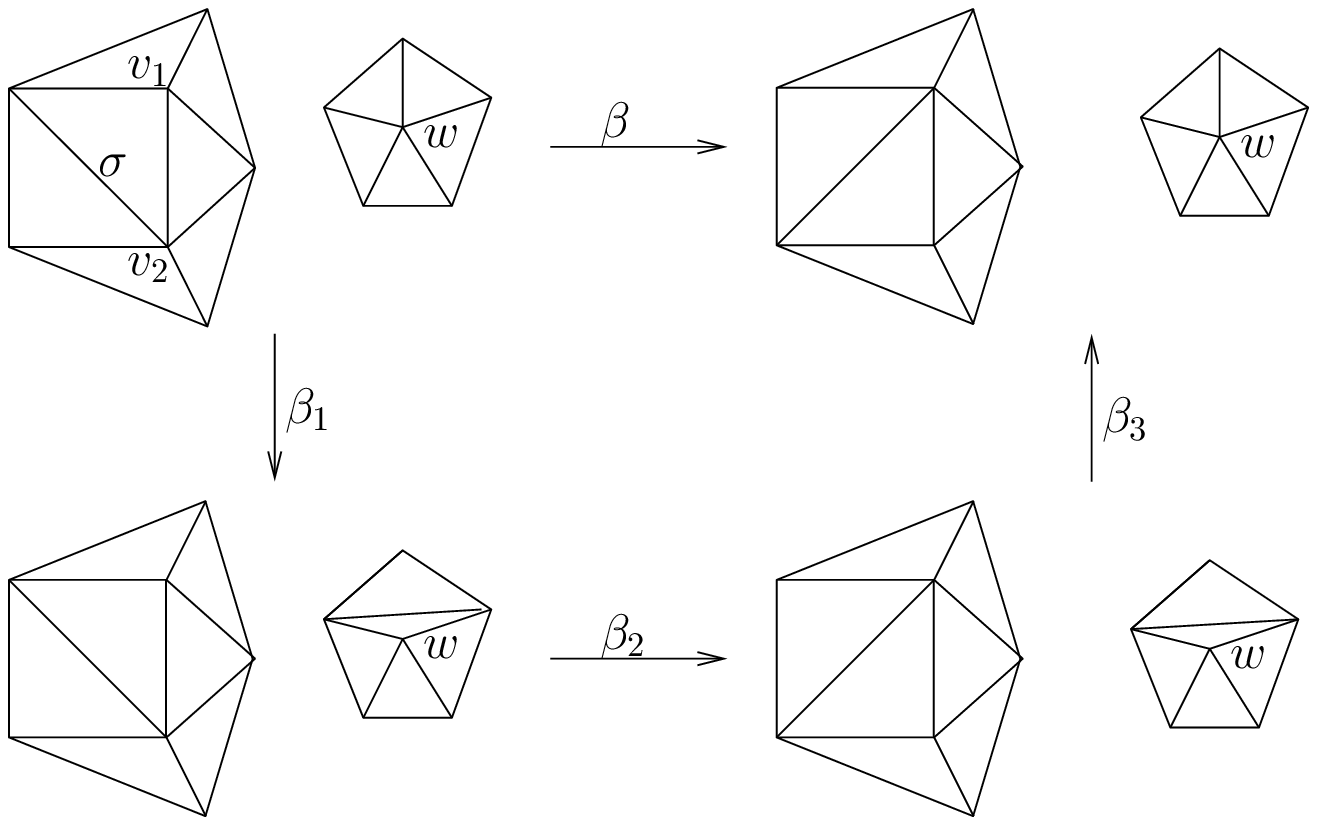}
					\caption{\textbf{2. There are two vertices $v_1$ and $v_2$, participating in~$\beta$ such that $\deg_{L_1}v_1 = 4$, $\deg_{L_2}v_2 = 4$.} Moreover, $L_1$ contains not more than two vertices of degree 4, as they should participate in the move $\beta$. Using Euler characteristic it is easy to show that in this case there are at least $8$ vertices of degree $5$ in $L_1$. Then at least one of these vertices does not belong to $\link v_1 \bigcup \link v_2$, denote such a vertex by $w$. Let us consider the pentagon $\link w$. $L_1$ can contain at most 2 of its diagonals. Then there is an edge adjacent to $w$ such that the move, associated with this edge, is defined. Hence, all the cycle obtained by commuting this move and the move $\beta$ is defined. So, we receive three moves $\beta_1$, $\beta_2$ and $\beta_3$, where each of them can be represented as a linear combination of elemenary cycles and bistellar flips with lower complexity according to the previous case. The original move $\beta$ can be represented in the same way.}
				\end{SCfigure}

			\subparagraph{$b=5$.}
				
				Each of the combinatorial spheres $L_1$ and $L_2$ does not contain vertices of degrees $3$ and $4$. In this case $L_1$ contains at least $12$ vertices of degree $5$. Among these 12 vertices there necessarily is a vertex $w$ such that it does not participate in $\beta$. Denote the vertices of $\link w$ by $u_1, u_2, u_3, u_4$ and $u_5$ in any cyclic order around $w$. From the five diagonals in the pentagon $\link w$ not more than two are present in each of the combinatorial spheres $L_1$ and $L_2$. This means that there is at least one diagonal, missing from both spheres. Without loss of generality, let this diagonal be $u_2u_5$. Then the elementary cycle $\gamma(L_1,\sigma,wu_1)$ is defined and the chain $\beta-\gamma(L_1,\sigma,wu_1)$ has its support in the graph $\Gamma_2^{k+\frac{2}{3}}$.

		\addtocounter{figure}{-4}
				
		\paragraph{The case of even $b$\,.}
					
			If $b$ is even bistellar moves with the biggest complexity split in pairs of successive moves:
			\[	
				\begin{diagram}
					\node[2]{L} \arrow{se,t}{\beta_2}
					\\
					\node[1]{L_1} \arrow{ne,t}{\beta_1}
					\node[2]{L_2}
				\end{diagram}
			\]
			Here we have $a(L) > a(L_1)$\,, $a(L) > a(L_2)$\,, $a(\beta_1)=a(\beta_2)=a(L)$.
				\\	
			Let $\beta_1=\beta^{-1}_{L,\sigma_1}$\,, $\beta_2=\beta_{L,\sigma_2}$\,.
			\subparagraph{$b=0.$}
				The moves $\beta_1^{-1}$ and $\beta_2$ reduce the number of vertices. The cycle $\gamma(L_1,\sigma_1,\sigma_2)$ is defined always, if $L$ has more than 5 vertices. If $L$ has 5 vertices, than the bistellar moves $\beta_1^{-1}$ and $\beta_2$ are equivalent.
		
			\subparagraph{$b=2.$}
				The complexes $L_1$ and $L_2$ contain vertices of degree $3$, but $L$ does not contain any vertices of degree $3$, and contains vertices of degree $4$. Then $\sigma_1$ and $\sigma_2$ are edges adjacent to some vertices $v_1$ and $v_2$ of degree {4}.
				
\begin{SCfigure}[5][hp] \label{extra2}
					\captionsetup{labelformat=empty}
					\includegraphics[scale=0.43]{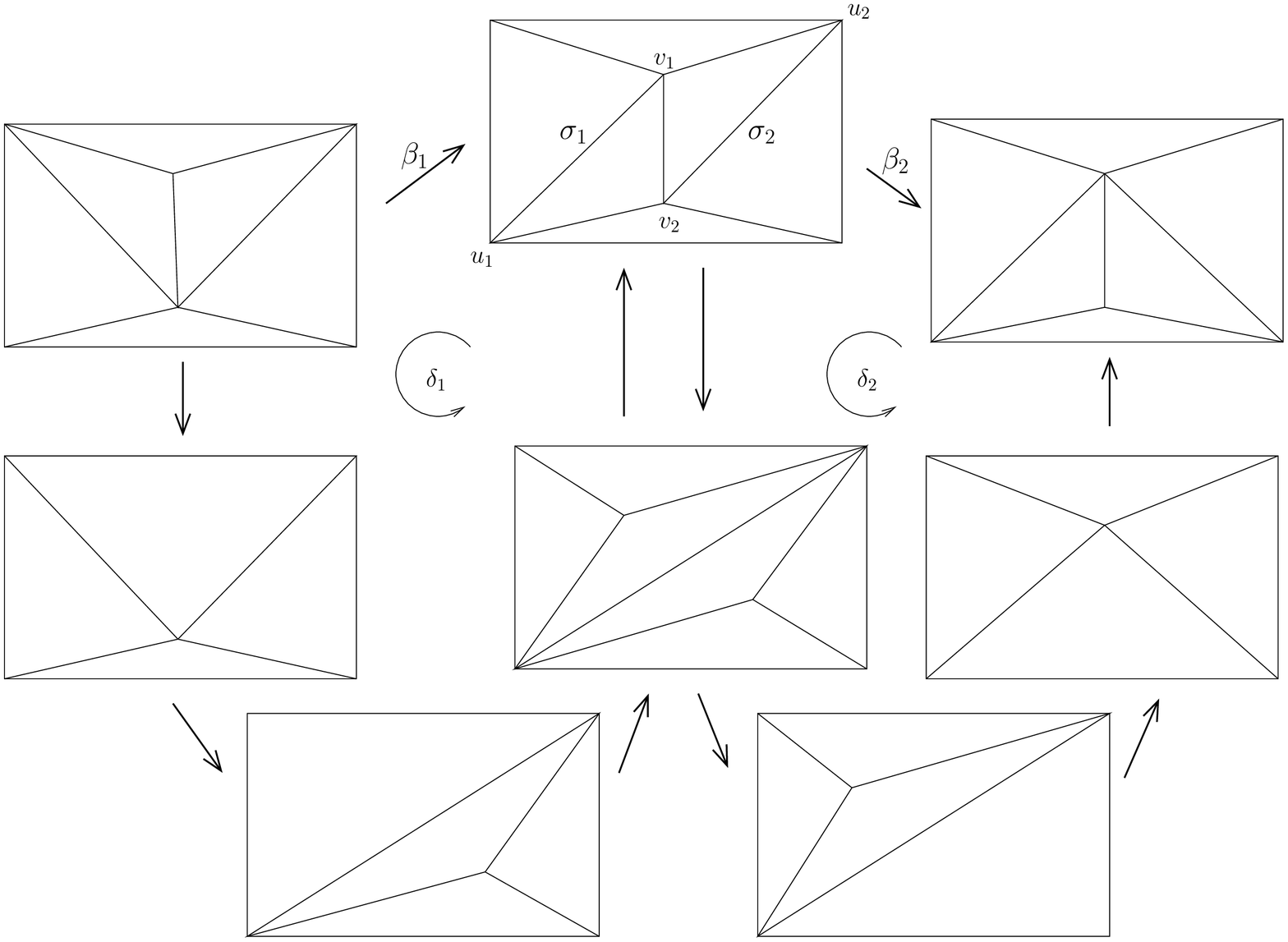}
					\caption{\textbf{4.} There is a unique case when the subtraction of the cycle $\gamma(L,\sigma_1,\sigma_2)$ from the chain $\beta_1+\beta_2$ does not lower the complexity of the chain. \textbf{This happens if the vertices $v_1$ and $v_2$ participate in both moves $\beta_1$ and $\beta_2$.} In this case the complexity of $\beta_1+\beta_2-\gamma(L_1,\sigma_1,\sigma_2)$ does not become lower than the complexity of the initial chain $\beta_1+\beta_2$. If the edge denoted $u_1u_2$ is not present in $L$, then two elementary moves of the second type $\delta_1$ and $\delta_2$ are defined. All new moves except for two cancelling ones have lower complexity than $a$, ie the support of the chain $\beta_1+\beta_2-\delta_1-\delta_2$ belongs to the graph $\Gamma_{2}^{k+\frac{1}{6}}$.}
				\end{SCfigure}

				If the cycle $\gamma(L,\sigma_1,\sigma_2)$ is defined, than $\beta_1 + \beta_2- \gamma(L,\sigma_1,\sigma_2) \in \Gamma_{2}^{k-\frac{1}{6}}$ except one case \ref{extra2} described last.
					
				The cycle $\gamma(L,\sigma_1,\sigma_2)$ is not defined in the following cases:
				
				\begin{SCfigure}[50][h]
					\captionsetup{labelformat=empty}
					\includegraphics[scale=0.45]{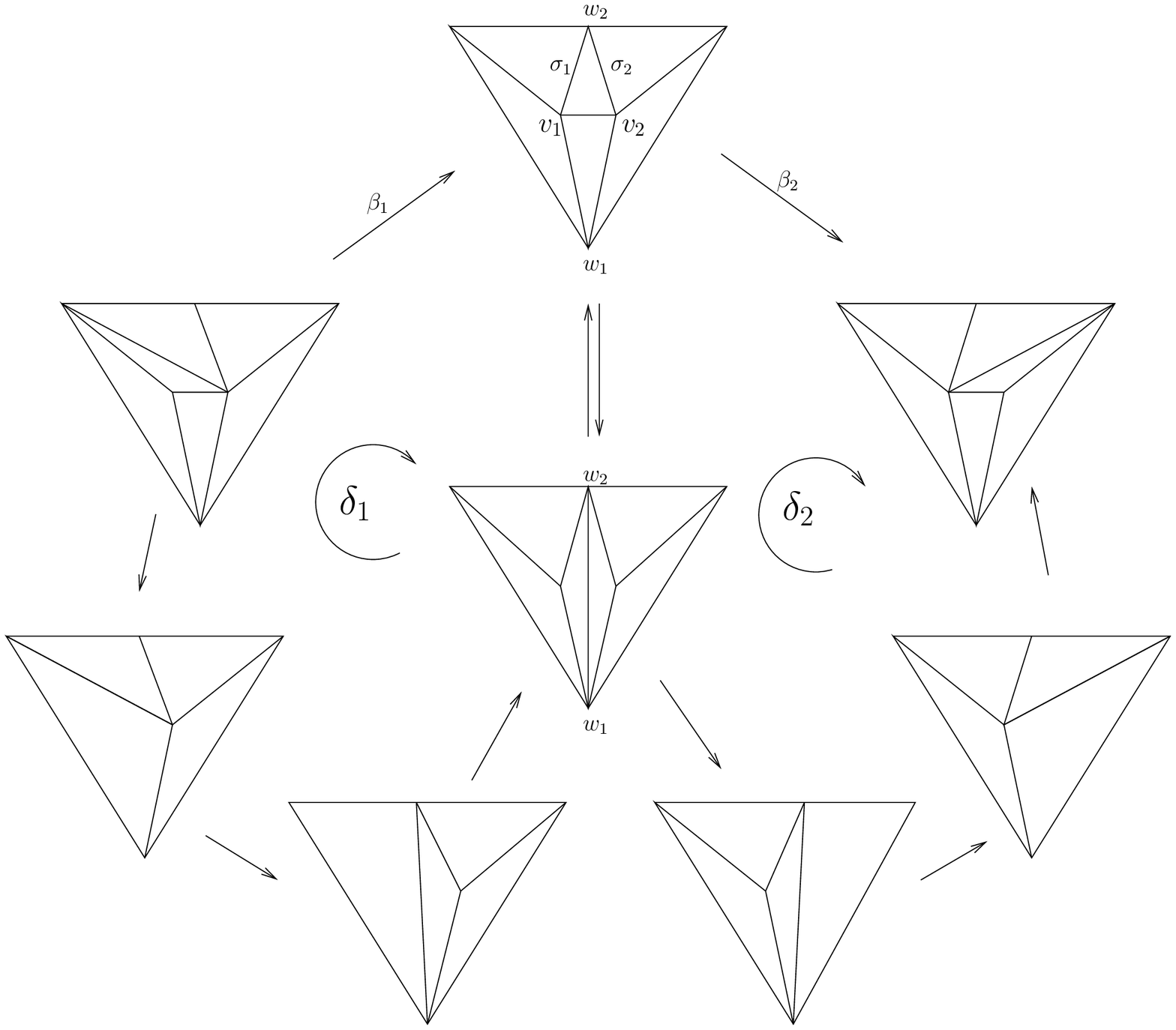}
					\caption{\textbf{1. The edges $\sigma_1$ and $\sigma_2$ are contained in a common triangle of the combinatorial sphere $L$, and their common vertex has a degree exceeding 4.} If the combinatorial there $L$ does not contain the edge $w_1w_2$ (the notations are on the picture), then we can apply the composition of two elementary cycles of the second type $\delta_1$ and $\delta_2$. Then new bistellar moves, except two cancelling out, will be less complex, and the support of the difference $\beta_1 + \beta_2 - \delta_1 - \delta_2$ lies in the graph $\Gamma_{2}^{k+\frac{1}{6}}$.}
				\end{SCfigure}
				
				\begin{SCfigure}[5][h]
					\captionsetup{labelformat=empty}
					\includegraphics[scale=0.46]{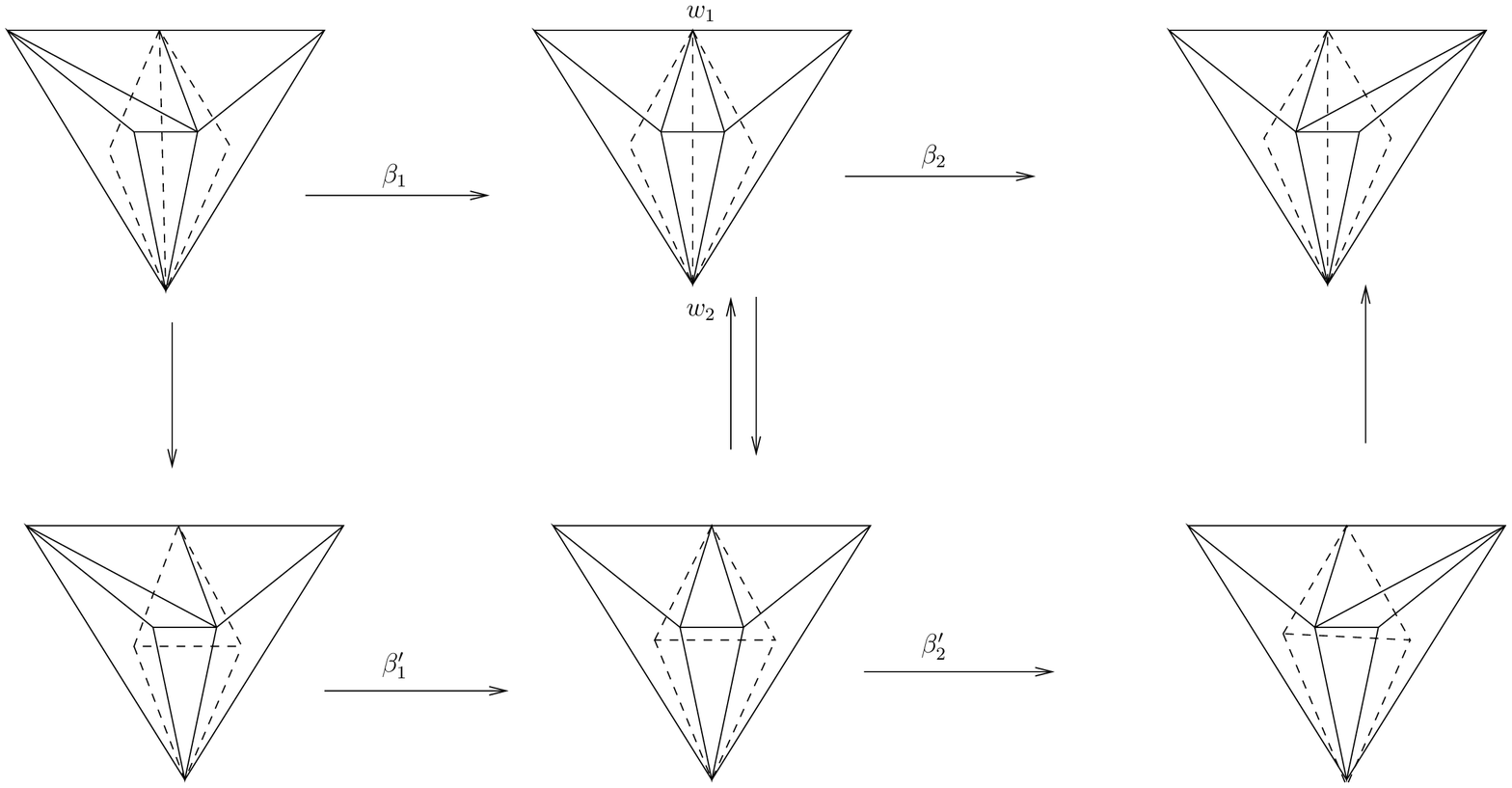}
					\caption{If the edge $w_1w_2$ already lies in $L$, then we can delete it, applying a composition of two elementary cycles of the first type, as shown on the figure. New vertical bistellar moves, except two cancelling out, have a complexity less than $a$, and new horizontal moves -- not exceeding $a$. The chain $\beta_1 + \beta_2 + \gamma(L,\sigma_1,w_1,w_2) - \gamma(L,\sigma_2,w_1,w_2)$ can be represented in the desired form using the previous case for the sum $\beta'_1 + \beta'_2$. Then $\beta_1 + \beta_2$ is represented as the sum of the result for $\beta'_1 + \beta'_2$, two elementary cycles and vertical moves with smaller complexities.}
				\end{SCfigure}
				
				\begin{SCfigure}[50][hp]
					\captionsetup{labelformat=empty}
					\includegraphics[scale=0.5]{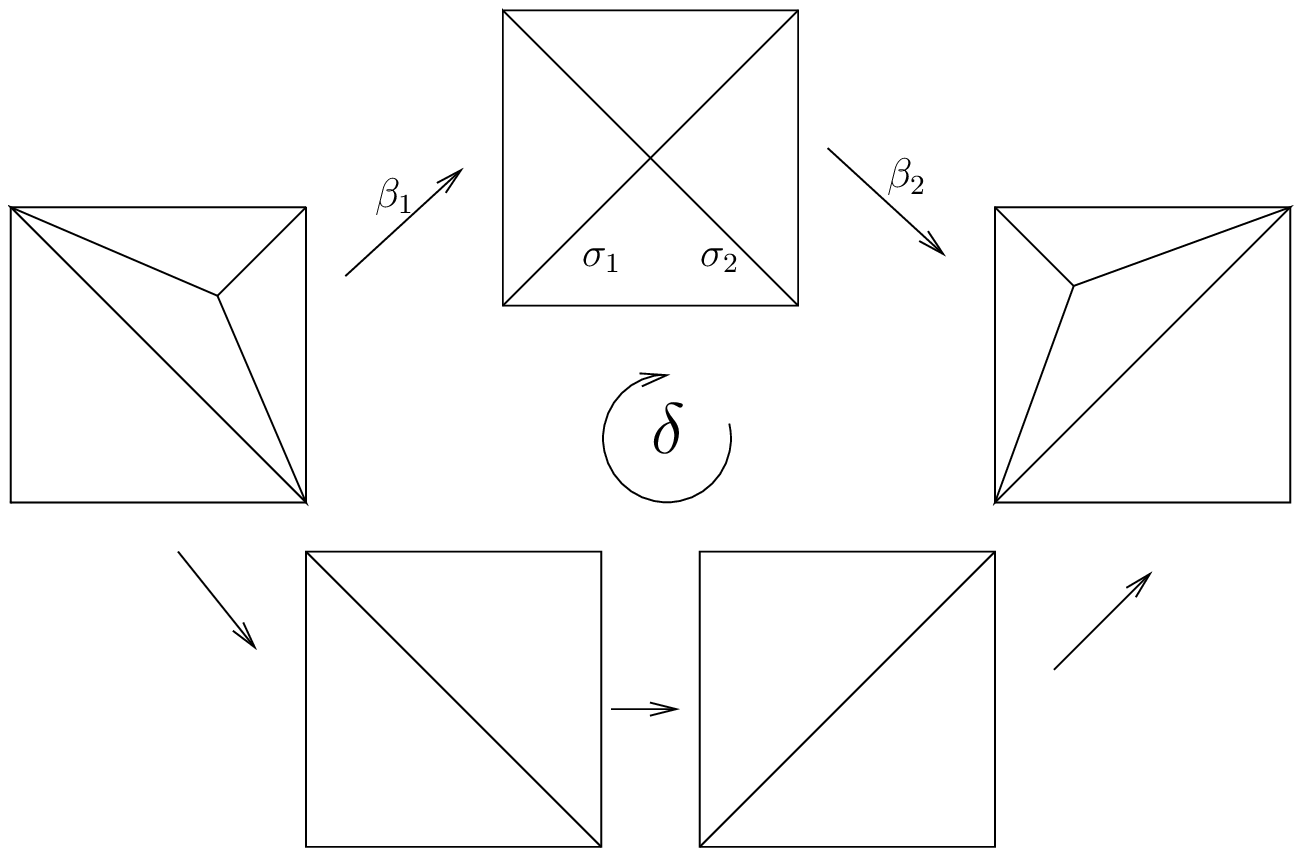}
					\caption{\textbf{2. The edges $\sigma_1$ and $\sigma_2$ are contained in a common simplex of the combinatorial sphere $L$, and the mutual vertex of these edges is of degree $4$.} In this case the elementary cycle $\delta$ of the second type is defined, and the support of the chain $\beta_1+\beta_2-\delta$ belongs to the graph $\Gamma_{2}^{k+\frac{1}{6}}$.}
				\end{SCfigure}
				
				\begin{SCfigure}[50][hp]
					\captionsetup{labelformat=empty}
					\includegraphics[scale=0.5]{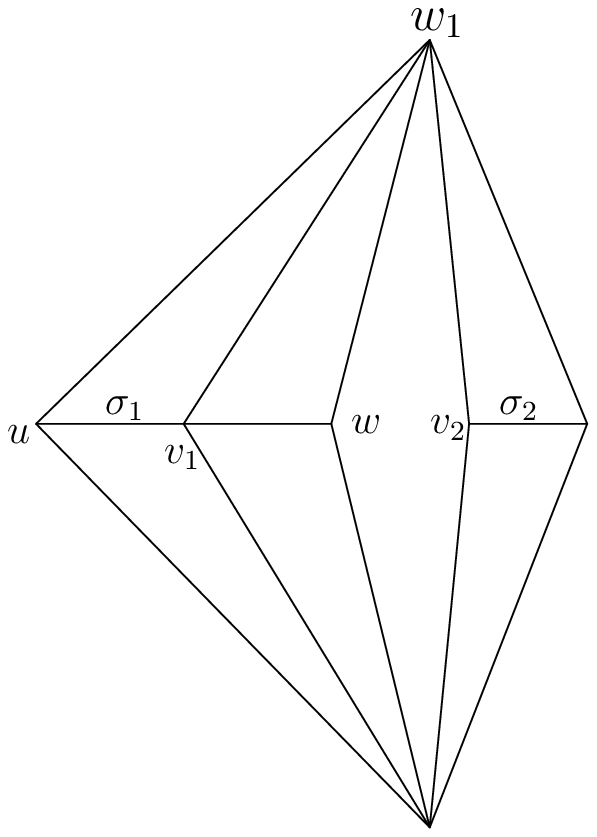}
					\caption{\textbf{3. The edges $\sigma_1$ and $\sigma_2$ are not contained in any common simplex, but their links in $L$ coincide, so the commutation of the bistellar moves $\beta_{\sigma_1}$ and $\beta_{\sigma_2}$ is impossible.} Then there is two different possibilities: $\sigma_1$ and $\sigma_2$ can have or not a common vertex. Suppose that these edges do not intersect. Then let $w$ be the vertex as in the figure (it is possible that $w=v_2$, this does not change the step of the algorithm). The edge $uw$ can not belong to the combinatorial sphere $L$, thus the move $\beta_3 = \beta_{v_1w_1}$ is defined. According to the previous case, the chain $\beta_1+\beta_3$ can be decomposed in a linear combination of elementary cycles and a chain with its support belonging to the graph $\Gamma^{k+\frac{1}{6}}$. The difference $\beta_2-\beta_3$ can be decomposed using the cycle $\gamma(L,\sigma_2,v_1w_1)$. Then the chain $\beta_1+\beta_2$ can be decomposed in the same way as $\beta_1+\beta_2 = (\beta_1+\beta_3) + (\beta_2-\beta_3)$}
				\end{SCfigure}
				
				\begin{SCfigure}[50][hp]
					\captionsetup{labelformat=empty}
					\includegraphics[scale=0.65]{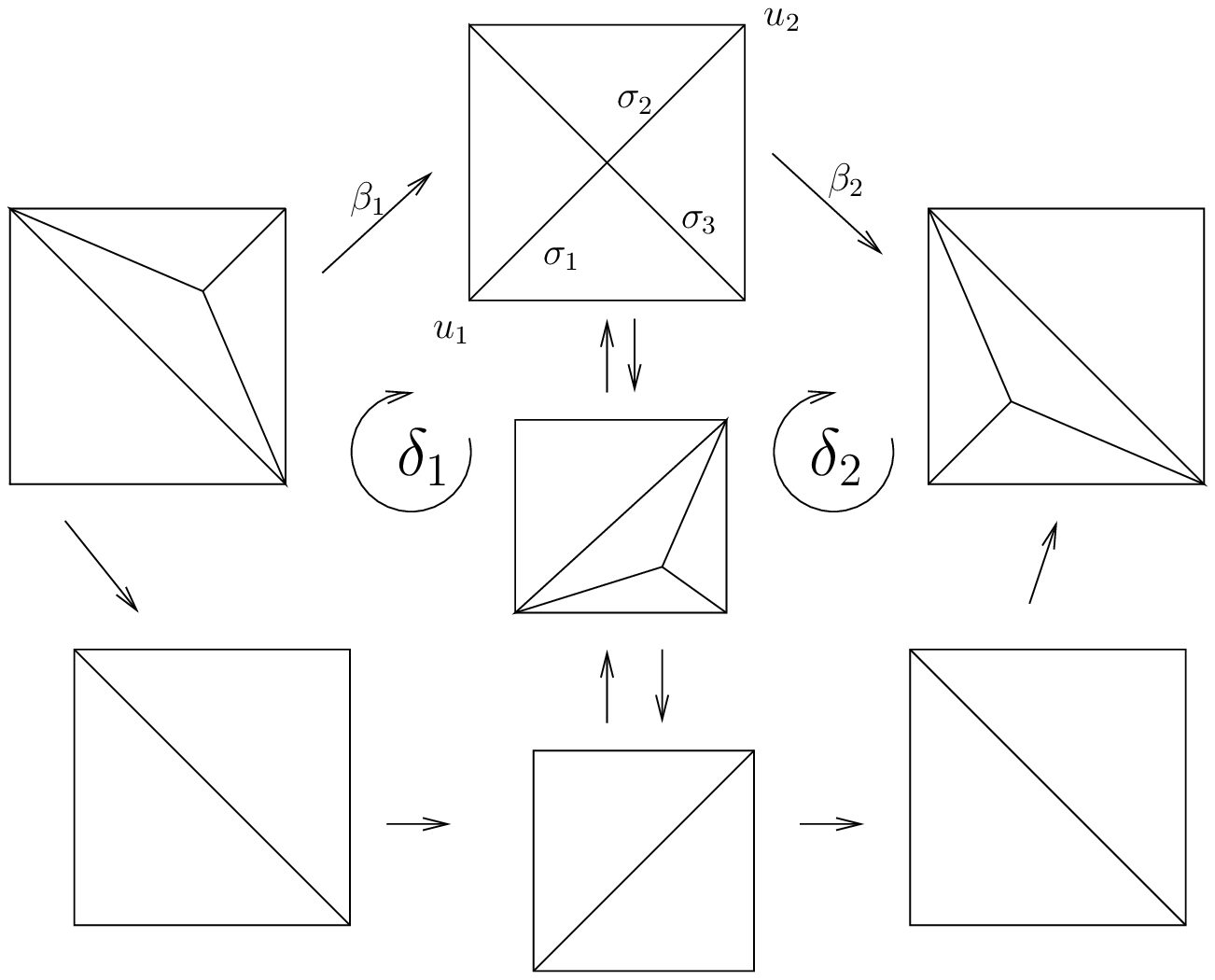}
					\caption{Now consider the case when $\sigma_1$ and $\sigma_2$ have a common vertex. If the diagonal $u_1u_2$ of the depicted quadrangle does not belong to $L$ then two elementary moves of the second type $\delta_1$ and $\delta_2$ are defined, and the complexity of all new moves except the cancelling ones is lower than the complexity of $\beta_1$ and $\beta_2$, hence the support of $\beta_1+\beta_2 - \delta_1 - \delta_2$ belongs to the graph $\Gamma_{2}^{k+\frac{1}{6}}$.}
				\end{SCfigure}
				
				\begin{figure}[hp]
					\captionsetup{labelformat=empty}
					\includegraphics[scale=0.5]{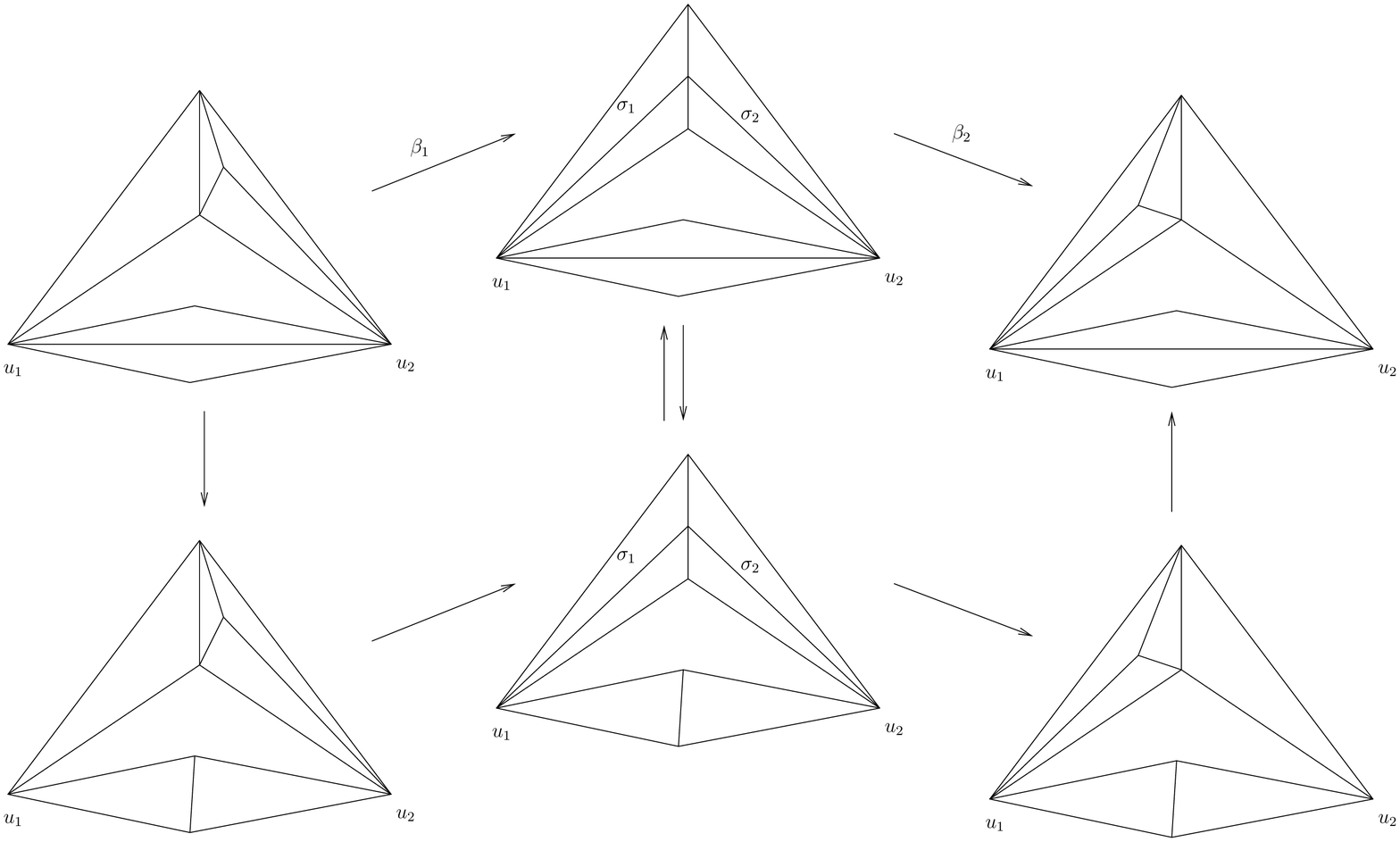}
					\caption{Now suppose that the diagonal $u_1u_2$ is present in the sphere $L$. This case can be solved in the same way as in case (1) (see figure above). Elementary cycles $\gamma(L,\sigma_1,w_1w_2)$ and  $\gamma(L,\sigma_2,w_1w_2)$ are defined. The chain $\beta_1 + \beta_2 + \gamma(L_1,\sigma_1,w_1w_2) - \gamma(L,\sigma_2,w_1w_2)$ can then be represented as a sum of moves with complexity less than $a$ and two moves that can be represented in the desired way according to the precious case. So, we described all the cases when the cycle $\gamma(L_1,\sigma_1,\sigma_2)$ is not defined.}
				\end{figure}

				\begin{SCfigure}[5][h!] 
					\captionsetup{labelformat=empty}
					\includegraphics[scale=0.43]{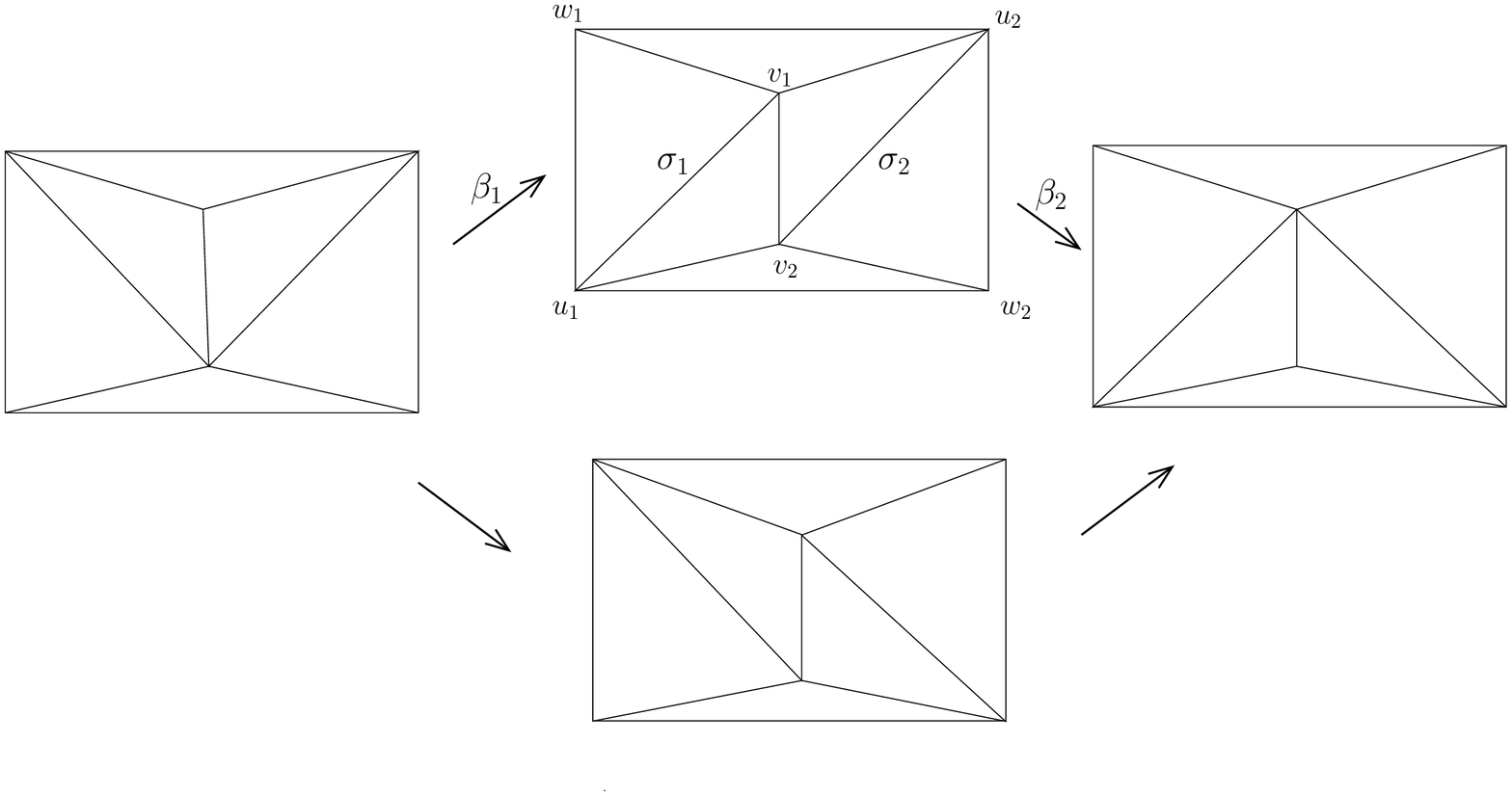}
					\caption{If the edge $u_1u_2$ is present in the combinatorial sphere $L$, then $L$ can not contain the edge denoted by $w_1w_2$. Consider the chain $\beta_1+\beta_2-\gamma(L_1,\sigma_1,\sigma_2)$. It can be represented in the desired way, as $w_1w_2$ does not belong to $L$ and we can use the previous case. Then $\beta_1+\beta_2$ can be represented as a linear combination of elementary moves and moves with lower complexity.}
				\end{SCfigure}
			
			\newpage
			
			\subparagraph{$\mathbf{b=4}$.}
				
				The combinatorial sphere $L$ does not contain any vertices of degree 3 and 4, one of the vertices of both $\sigma_1$ and $\sigma_2$ is of degree $5$. Let $v_1$ and $v_2$, respectively, be those vertices. $L$ contains not less than $12$ vertices of degree $5$, wherein not more than 8 vertices participate in the moves $\beta_1$ and $\beta_2$. Hence there is a common vertex of degree $5$ in the three combinatorial spheres $L_1$, $L_2$ and $L_3$. Denote this vertex by $v$. Among the $5$ edges adjacent to $v$ there are at least $3$ edges such that moves associated with these edges are defined. Denote these vertices by $e_1$, $e_2$ and $e_3$. Then the cycle $\gamma(L_1,\sigma_1,e_i)$, as well as the cycle $\gamma(L,\sigma_2,e_i)$, is defined for at least two of three edges $e_i$. Thus there is an $i$ such that both cycles $\gamma(L_1,\sigma_1,e_i)$ and $\gamma(L,\sigma_2,e_i)$ are defined. Then the support of the chain $\beta_1 + \beta_2 - \gamma(L_1,\sigma_1,e_i) + \gamma(L,\sigma_2,e_i)$ belongs to the graph $\Gamma^{k+\frac{1}{2}}$ except one case, similar to the case (4) for $b=2$.
				
				\begin{SCfigure}[5][!hb]
					\captionsetup{labelformat=empty} 
					\includegraphics[scale=0.6]{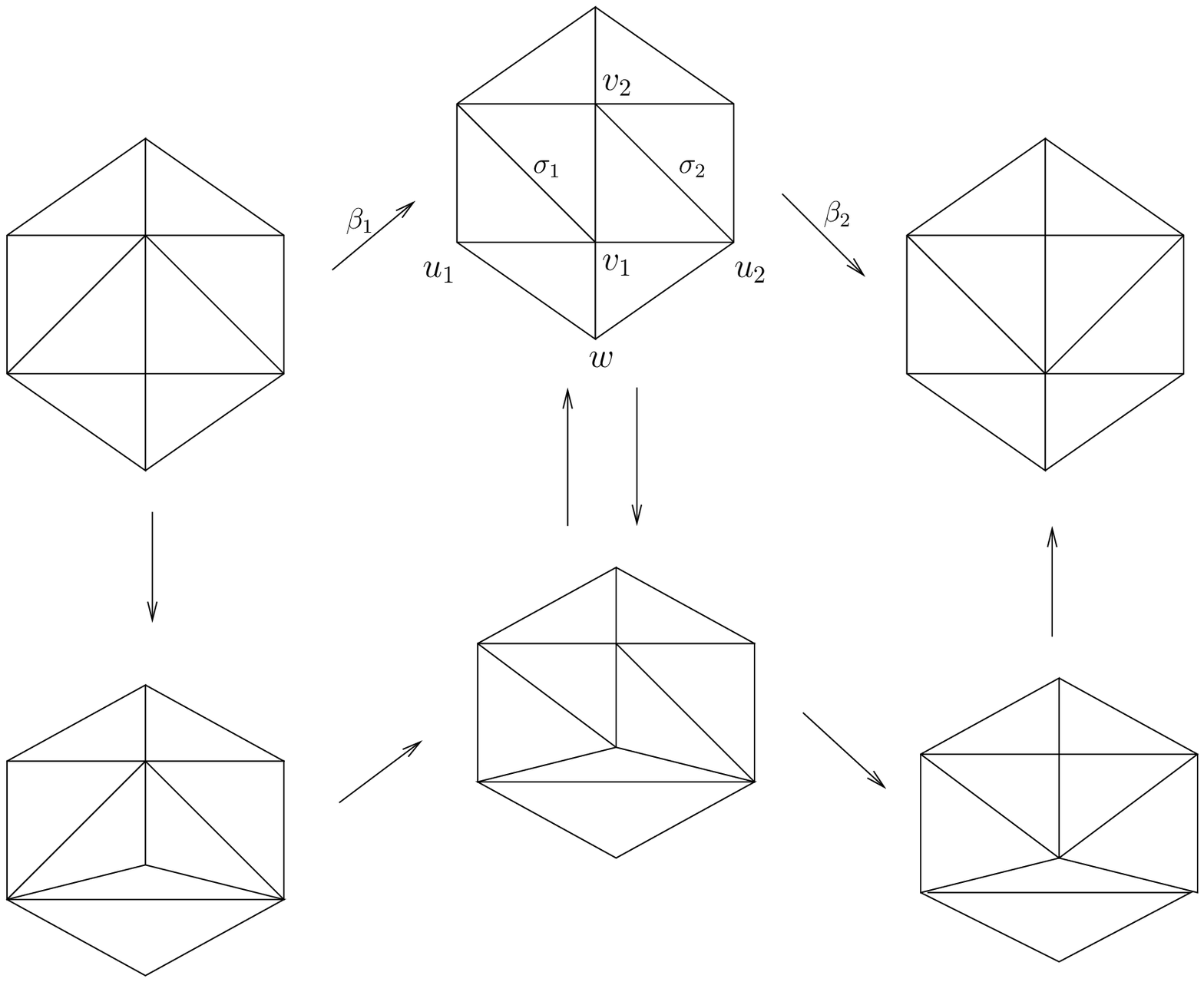}
					\caption{The last case appears if the degrees of the vertices $v_1$ and $v_2$ do not decrease to $4$ under the moves of the chain $\beta_1 + \beta_2 - \gamma(L_1,\sigma_1,e_i) + \gamma(L,\sigma_2,e_i)$. Then the vertices $v_1$ and $v_2$ belong to a common edge, as on the figure. If $L$ does not contain the edge denoted by $u_1u_2$, then elementary cycles of the first type $\gamma(L, \sigma_1, v_1w)$ and $\gamma(L, \sigma_2, v_1w)$ are defined, and the support of the chain $\beta_1 + \beta_2- \gamma(L,\sigma_1,v_1w) + \gamma(L,\sigma_2,v_1w)$ belongs to the graph $\Gamma^{k+\frac{1}{2}}$. If $L$ contains the edge $u_1u_2$, then, as in the case (3) for $b=2$, the cycles $\gamma(L,\sigma_1,u_1u_2)$ and $\gamma(L,\sigma_2,u_1u_2)$ are defined and the chain $\beta_1+\beta_2+\gamma(L,\sigma_1,u_1u_2)-\gamma(L,\sigma_2,u_1u_2)$ is represented as a sum of moves with lower complexities and two moves, where the edge $u_1u_2$ is absent.}
				\end{SCfigure}
						
				We proved the theorem stating that
				\begin{theorem}
				Any cycle in the graph $\Gamma_2$ can be represented as a linear combination of elementary cycles.
				\end{theorem}
				
				This theorem has also been proved by Gaifullin \cite{G04} using Steinitz theorem, but the proof here is necessary for the realization as it is completely explicit. The subcases for $b=2$ and $b=4$ where the cycle $\beta_1 + \beta_2 - \gamma(L_1,\sigma_1,e_i) + \gamma(L,\sigma_2,e_i)$ is not defined, as well as one subcase for $b=3$ were added to complete the algorithm from \cite{G08}.
				
				\addtocounter{figure}{-9}
	
\section{The realization of the algorithm}
\label{program}
	
	In the previous sections the proof of Theorem \ref{main} has been reduced to the computation of $p_1(M^8_{15})$. We will do this using the described algorithm.

	$M^8_{15}$ has $3003$\, $4$-simplexes. Though some of these simplices can be taken to each other by automorphisms of $M^8_{15}$, there still will be more than 60 combinatorial types of $\link \sigma^4$. Hence the computation by hand is labor intensive. But as this algorithm is completely combinatorial, it can be realized on a computer.
	
	Checking if two given combinatorial spheres are isomorphic is a computationally hard problem. Gaifullin's algorithm operates with isomorphism classes of combinatorial manifolds (and bistellar moves). We would like to avoid checks of sphere isomorphism for the program to work faster and the realization to be easier. Let us introduce an additional construction for this purpose based on the graph $\Gamma_2$.  Define the graph $\widetilde{\Gamma}_2$ as follows. This graph has as vertices oriented combinatorial 2-spheres with vertices labeled by pairwise distinct natural numbers (not necessarily successive), up to label preserving isomorphism, and its edges are equivalence classes of bistellar moves, preserving orientation and respecting the labeling of vertices. If a vertex is added under a bistellar move then it can have any possible label.
	
\begin{figure}[h]
		\begin{center}
		\includegraphics[scale=0.5]{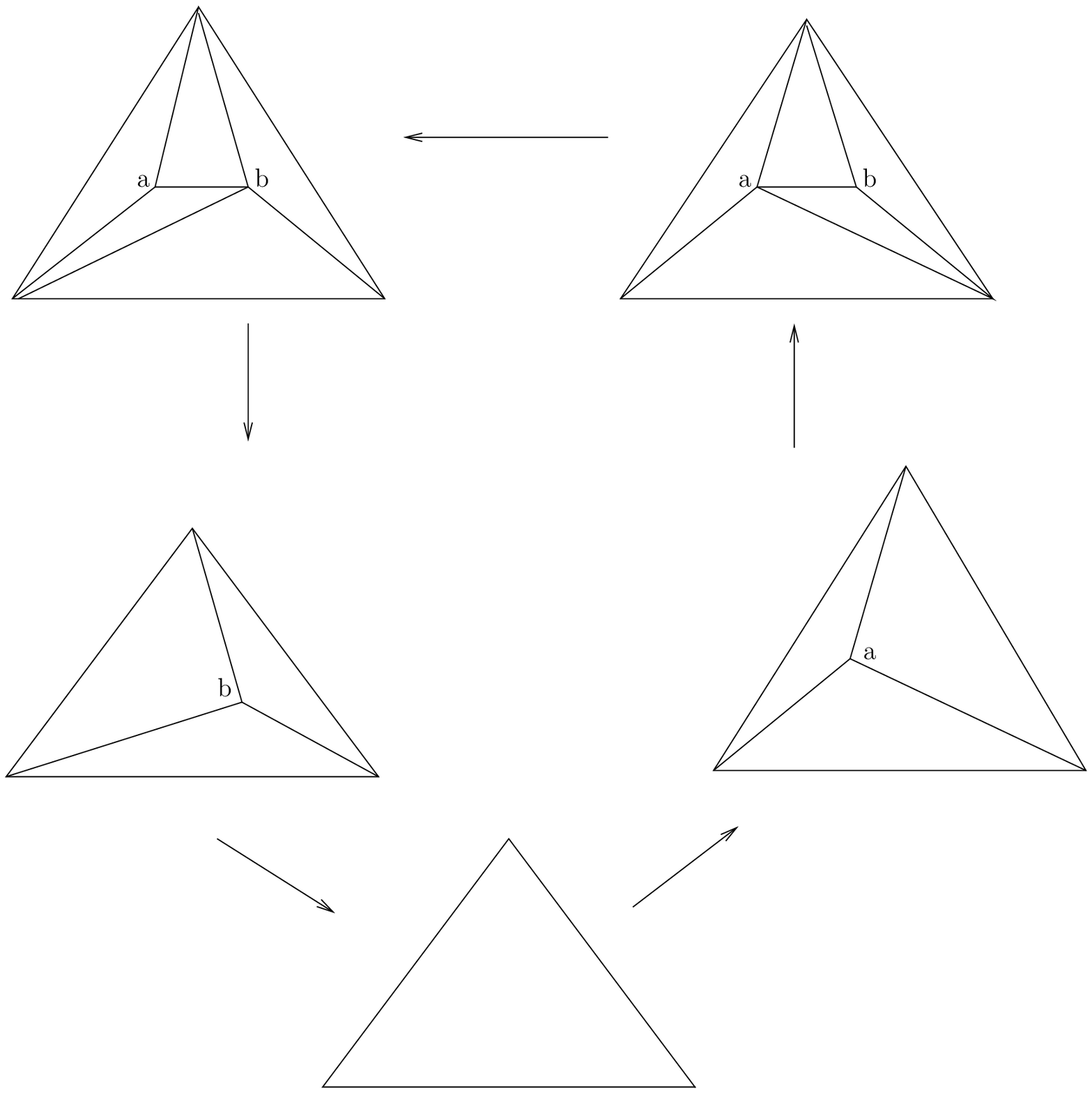}
		\end{center}
		\caption{The changed elementary cycle of type (2a)}
		\label{2a}
	\end{figure}

	There is a natural map $p \colon \widetilde{\Gamma}_2 \longrightarrow \Gamma_2$, that forgets the vertex labeling of the sphere. The pull-back $p^* \colon C^1(\Gamma_2,\mathbb{Q}) \longrightarrow C^1(\widetilde{\Gamma}_2,\mathbb{Q})$ sends the cocycle $h$ to a cocyle $\tilde{h}$. We will call elementary cycles in $\widetilde{\Gamma}_2$ the same cycles that were elementary in $\Gamma_2$, but with a fixed labeling of the vertices of all combinatorial spheres, such that every move is well-defined as an edge in $\widetilde{\Gamma}_2$. The only exception will be the cycle (2a), as it is impossible to label vertices in the figure in a way for all moves to respect the labeling. We shall add to this cycle two inverse bistellar moves as on Fig.\ref{2a} for this cycle to be defined on $\widetilde{\Gamma}_2$.
	
	The image of any elementary cycle in $\widetilde{\Gamma}_2$ under the map $p$ is an elementary cycle in $\Gamma_2$. The algorithm of cycle decomposition in the graph $\Gamma_2$ is naturally transorted on the graph $\widetilde{\Gamma}_2$. As $\tilde{h}=p^*(h)$, the value of $\tilde{h}$ on an elementary cycle in $\widetilde{\Gamma}_2$ is equal to the value of $h$ on the image of this elementary cycle under the projection on $\Gamma_2$. Hence the values of the cocycle $\tilde{h}$ on elementary cycles are computed in the same way as the values of $h$.
	
	\newpage

	Let the vertices of the initial complex $K$ be labeled. Consider now the steps of the realization of the first Pontryagin class computation algorithm.
	The algorithm for a labeled complex consists of the following steps:
		\begin{enumerate}
			\item For every oriented $(n-4)$--simplex $\sigma$ of the complex $K$ find a sequence of moves $\xi_\sigma$, respecting the labeling, that transform the link (with induced orientation) of $\sigma$ into the boundary of a simplex. This step is realized with the help of the program \texttt{BISTELLAR} \cite{BISTELLAR} (the programming language is \texttt{GAP} \cite{GAP4}). The algorithm used in this program is not a full algorithm checking the isomorphism of a given complex and a combinatorial sphere (even in the 3-dimensional case) because no estimations are known on the time of work of the program, but it works effectively in all arising examples. The program \texttt{BISTELLAR} explicitly finds a sequence of bistellar moves, gradually decreasing the number of vertices of the complex. In the case of a combinatorial sphere this allows to descend to the boundary of a simplex.
							
			\item For each vertex $v$ of the link of every $(n-4)$-simplex $\sigma$ (as well as all new vertices appearing in $\xi_\sigma$) consider $\link_{\link\sigma}(v)$. Then each of the obtained complexes is a combinatorial 2-sphere. Induce the sequences $\xi_\sigma$ of bistellar moves on these complexes as on subcomplexes of $\link\sigma$ preserving the labeling. Denote the sequence induced on the subcomplex $\link_{\link\sigma}(v)$ by $\xi_{\sigma,v}$. We should be careful about vertices that can be added to $\link_\sigma$ in moves used in $\xi_\sigma$, these new vertices shall also be considered. Denote by $V(\xi_\sigma)$ the set of all vertices that appear in the moves of the chain $\xi_\sigma$.

			\item Let us choose a natural way to construct a chain $\kappa(L)$ of moves between a combinatorial $2$-sphere $L$ and the boundary of the $3$-simplex (ie if two isomorphic combinatorial spheres $L_1$ and $L_2$ have the same labelings, then the chosen chains will be isomorphic and identically labeled). For example, we can apply the lexicographically first possible bistellar move decreasing the complexity of the combinatorial sphere $L$ and in the same way descend to the boundary of the simplex. For each chain of moves $\xi_{\sigma,v}$, reducing a combinatorial 2-sphere to the boundary of a 3-simplex, we have the chain $\xi_{\sigma,v}-\kappa(\link_{\link\sigma}(v))$ from $\partial \Delta^3$ to $\partial \Delta^3$. The resulting simplexes can be labeled in different ways. But boundaries of the simplex $\partial \Delta^3$ labeled in different ways can be joined with a sequence of moves respecting orientation in the following way. Denote our boundaries of the simplex by $\partial \Delta_1^3$ (with labels $u_1$, $v_1$, $w_1$ and $z_1$) and $\partial \Delta_2^3$ (with labels $u_2$, $v_2$, $w_2$ and $z_2$). We will change the labels of the vertices one by one, for example, let us show the sequence that changes the label $u_1$ into $u_2$. This will be a sequence consisting of three moves (add a vertex with label $u_2$, then make the vertex labeled $u_1$ be of degree $3$, then remove it): 
			\begin{align*}
			& \{\{u_1,v_1,w_1\},\{u_1,v_1,z_1\},\{u_1,w_1,z_1\},\{v_1,w_1,z_1\}\}\longrightarrow \\ \longrightarrow & \{\{u_1,w_1,u_2\},\{u_1,v_1,u_2\},\{v_1,w_1,u_2\},\{u_1,v_1,z_1\},\{u_1,w_1,z_1\},\{v_1,w_1,z_1\}\} \longrightarrow \\ \longrightarrow & \{\{u_1,w_1,u_2\},\{u_1,z_1,u_2\},\{v_1,w_1,u_2\},\{u_2,v_1,z_1\},\{u_1,w_1,z_1\},\{v_1,w_1,z_1\}\} \longrightarrow \\ \longrightarrow &
			\{\{u_2,v_1,w_1\},\{u_2,v_1,z_1\},\{u_2,w_1,z_1\},\{v_1,w_1,z_1\}\}
			\end{align*}
			Moreover, we did not use combinatorial spheres with more than $5$ vertices. It is easy to verify that new moves constructing the chain between differently numerated $\partial \Delta^3$ give no contribution the value of the formula. Denote the chain joining $\partial \Delta_1^3$ and $\partial \Delta_2^3$ by $\zeta(\Delta_1,\Delta_2)$ Then we have a cycle in $\widetilde{\Gamma}_2$ $$\xi_{\sigma,v}-\kappa(\link_{\link\sigma}(v)) + \zeta(\Delta_1,\Delta_2) \in Z^1(\Gamma_2,\mathbb{Q})$$. Denote this cycle by $\eta_{\sigma,v}$.
								
			\item The resulting cycles $\eta_{\sigma,v}$ are cycles in the graph $\widetilde{\Gamma}_2$. Decompose them in a linear combination of elementary cycles.
								
			\item Compute the contribution of each elementary cycle. For each $\sigma$ receive its contribution $$f\left( \left\langle link\,\sigma\right\rangle \right) = \sum_{v \in V(\xi_\sigma)} \langle c,\eta_{\sigma,v}\rangle$$ and construct the cycle $$f_\sharp(L) = \sum_{\sigma \in L, \dim \sigma = n-4} f\left( \left\langle link\,\sigma\right\rangle \right)\, \sigma,$$ representing the homology element, dual to the first Pontryagin class. To receive more explicit results in cohomology groups and compute the first Pontryagin number we use the package \texttt{simpcomp}\cite{simpcomp}.
						
		\end{enumerate}
		
	We need a remark for the step 3. The problem of the constructed chain $\kappa$ is that it does not preserve the localness of the formula as it depends on the numeration of the complex
	
	\begin{lemma}[Gaifullin, \cite{G04}]
	
		The homology class in the computation of the first Pontryagin class does not depend on the choice of closure of the chain in the graph $\Gamma_2$ if the closure depends uniquely on the labeling of the 2-sphere.
	
	\end{lemma}
	
	The author wrote a program using the programming language \texttt{GAP}\cite{GAP4} that realizes the algorithm. It takes a simplicial complex and gives the first Pontryagin class as well as the dual to it.
	
	Launching the program for $M_{15}^8$ gave the following answer: the first Pontryagin class is proportional to the image of one of two generators of $H^4(M_{15}^8,\mathbb{Z})$ under the natural inclusion with coefficient $2$. This proves the result announced in the beginning of the paper:
	
	\begin{ntheorem}
			The first rational Pontryagin class $p_1(M^8_{15})$ is equal to $2u$ where $u$ is the image of one of two generators of the group $H^4(M^8_{15},\mathbb{Z})\cong \mathbb{Z}$ under the natural embedding $H^4(M^8_{15},\mathbb{Z})\subset H^4(M^8_{15},\mathbb{Q})$.
	\end{ntheorem}
	
	Hence with Proposition \ref{prop} we have the following result
	\begin{corollary}
		$M^8_{15}$, $\widetilde{M}\vphantom{M}^8_{15}$ and  $\widetilde{\widetilde{M}}\vphantom{M}^8_{15}$ are PL homeomorphic to $\mathbb{H}P^2$ and are minimal triangulations of $\mathbb{H}P^2$.
	\end{corollary}
	
The author would like to thank his advisor Alexander A. Gaifullin for suggesting this interesting problem, for invaluable discussions, constant attention to this work and patience.

\end{document}